\documentclass[11pt,english]{article}
\usepackage[latin9]{inputenc}
\usepackage{geometry}
\usepackage{color}
\usepackage{babel}
\usepackage{array}
\usepackage{prettyref}
\usepackage{float}
\usepackage{booktabs}
\usepackage{multirow}
\usepackage{amsmath}
\usepackage{amsthm}
\usepackage{amssymb}
\usepackage{graphicx}
\usepackage{setspace}
\usepackage[authoryear]{natbib}
\usepackage[unicode=true,pdfusetitle,
bookmarks=true,bookmarksnumbered=true,bookmarksopen=false,
breaklinks=false,pdfborder={0 0 1},backref=false,colorlinks=true]
{hyperref}
\hypersetup{
	citecolor=black}
\makeatletter

\providecommand{\tabularnewline}{\\}
\floatstyle{ruled}
\newfloat{algorithm}{tbp}{loa}
\providecommand{\algorithmname}{Algorithm}
\floatname{algorithm}{\protect\algorithmname}

\theoremstyle{remark}
\newtheorem*{notation*}{\protect\notationname}
\theoremstyle{plain}
\newtheorem{theorem}{\protect\theoremname}
\theoremstyle{plain}
\newtheorem{assumption}{\protect\assumptionname}
\theoremstyle{plain}
\newtheorem{proposition}{\protect\propositionname}
\theoremstyle{plain}
\newtheorem{lemma}{\protect\lemmaname}
\theoremstyle{remark}
\newtheorem{remark}{\protect\remarkname}
\theoremstyle{definition}
\newtheorem{definition}{\protect\definitionname}
\theoremstyle{definition}

\newtheorem{corollary}{Corollary} 

\usepackage{xcolor}

\hypersetup{
	colorlinks,
	linkcolor={red!60!black},
	citecolor={blue!60!black},
	urlcolor={blue!80!black}
}

\allowdisplaybreaks

\usepackage{algorithm,algpseudocode}

\usepackage{etoolbox}

\makeatother

\providecommand{\assumptionname}{Assumption}
\providecommand{\definitionname}{Definition}
\providecommand{\examplename}{Example}
\providecommand{\lemmaname}{Lemma}
\providecommand{\notationname}{Notation}
\providecommand{\propositionname}{Proposition}
\providecommand{\remarkname}{Remark}
\providecommand{\theoremname}{Theorem}

\usepackage{natbib}
\bibpunct[, ]{(}{)}{,}{a}{}{,}%
\usepackage{authblk}
\title{An Accelerated Fitted Value Iteration Algorithm for MDPs with Finite and Vector-Valued Action Space}
\author[1]{Sixiang Zhao}
\author[2]{William B. Haskell}
\author[3]{Michel-Alexandre Cardin}
\affil[1]{Sino-US Global Logistics Institute, Shanghai Jiao Tong University, Shanghai 200030, China, sixiang.zhao@sjtu.edu.cn}
\affil[2]{Krannert School of Management, Purdue University, West Lafayette, IN 47907, USA, whaskell@purdue.edu}
\affil[3]{Imperial College London, Dyson School of Design Engineering, London SW72AZ, UK, m.cardin@imperial.ac.uk}
\date{}                     
\begin{document}

	\maketitle
	\begin{abstract}
		This paper studies an accelerated fitted value iteration (FVI) algorithm to solve high-dimensional Markov decision processes (MDPs). FVI is an approximate dynamic programming algorithm that has desirable theoretical properties. However, it can be intractable when the action space is finite but vector-valued. To solve such MDPs via FVI, we first approximate the value functions by a two-layer neural network (NN) with rectified linear units (ReLU) being activation functions. We then verify that such approximators are strong enough for the MDP. To speed up the FVI, we recast the action selection problem as a two-stage stochastic programming problem, where the resulting recourse function comes from the two-layer NN. Then, the action selection problem is solved with a specialized multi-cut decomposition algorithm. More specifically, we design valid cuts by exploiting the structure of the approximated value functions to update the actions. We prove that the decomposition can find the global optimal solution in a finite number of iterations and the overall accelerated FVI is consistent. Finally, we verify the performance of the FVI algorithm via a multi-facility capacity investment problem (MCIP). A comprehensive numerical study is implemented, where the results show that the FVI is significantly accelerated without sacrificing too much in precision.
		
	\end{abstract}

	
	%
	
	
	\section{Introduction}
	
	Sequential decision problems under uncertainty are common in the community of production and operations management. In these problems, decision makers observe uncertain parameters and make decisions alternately and sequentially. Such problems can be found in inventory control \citep{feng_replenishment_2015}, capacity investment \citep{zhao_decision_2018}, patient scheduling \citep{yu_appointment_2020}, revenue management \citep{luo_revenue_2016}, to name a few. 
	
	Markov decision process (MDP) is one prevalent method to model sequential decision problems, and the MDP models, with either infinite or finite horizon, can be solved via value iteration (VI). VI is a dynamic programming (DP)-based algorithm which uses lookup table representations of the value functions and solve it in a recursive manner. However, VI becomes intractable when the state space of the MDP is multi-dimensional or continuous. To address this issue, one can use approximation. 
	
	Fitted value iteration (FVI), or approximate value iteration in some literature, is developed from the traditional VI. In FVI, the value functions are firstly approximated via some approximators, such as linear/piecewise linear functions, polynomial functions, or neural networks (NNs). The algorithm then proceeds via solving the Bellman equations recursively using the approximated value functions, instead of its exact representations. This type of algorithms is also known as approximate dynamic programming (ADP). \cite{munos_finite-time_2008} has provided desirable theoretical supports for FVI by showing when this algorithm will yield good policy. More specifically, if the approximators are rich enough and the MDP models satisfy some mild conditions, we can achieve $\epsilon$-optimal policy with high probability. 
	
	Though the FVI algorithm has desirable theoretical properties, it still has limitations in practice. That is, it assumes that the underlying MDP has a finite action space so that we can enumerate all possible actions. For a sample-based FVI, we need to solve the action selection problems for each state sample in each period once the value functions are approximated. This, however, can be extensively time-consuming when the action space is finite but vector-valued. 
	
	Our main contribution in this paper addresses the tension between the theory and solvability of FVI in the context of MDPs with finite and vector-valued action space by developing a novel acceleration scheme for the action selection procedure. Our acceleration algorithm for the action selection procedure requires that the transition function of the MDP model is linear with respect to the actions. This type of models are common in practice; inventory control problems, capacity expansion problems \citep{zhao_approximate_2017}, or dynamic assignment problems \citep{spivey_dynamic_2004} all fall into this community.
	
	We first propose a neural network-based FVI (abbreviated as NN-FVI) where the value functions are approximated via a two-layer NN with rectified linear units (ReLU) being activation functions. We verify that this choice of approximator is rich enough so that our FVI is consistent with some mild assumptions of the MDP. Then, we recast the action selection procedure of FVI as a two-stage stochastic programming model, where the recourse value function is returned from the NNs and so it may be non-concave. A two-stage stochastic programming with non-concave recourse function is notoriously difficult when we are maximizing the objective function. To address this, we further exploit the structure of the recourse function and propose a customized decomposition algorithm. In the algorithm, several types of valid inequalities are added to speed up the search of the optimal action. We prove that this acceleration scheme can converge to the optimal action with a finite number of iterations.
	
	To verify the performance of the proposed method, we apply it to a multi-facility capacity investment problem (MCIP). In particular, the capacity of the facilities can either be expanded or contracted, which is different from the exiting models where the capacity investment is irreversible \citep{zhao_approximate_2017, zhao_decision_2018}. We model this problem as an MDP and then analyze the properties of its value functions. We solve the model via the proposed FVI method and implement a comprehensive numerical study, where the results show that the FVI is significantly accelerated without sacrificing too much in precision.
	
	The remainder of the paper is organized as below. We review the related work in Section \ref{sec:Related_Work}. In Section \ref{sec:Preliminaries}, we present some preliminaries about MDP and the main assumptions in our paper. The proposed NN-FVI, along with the analysis of its consistency, are presented in Section \ref{sec:NNFVI}. Section \ref{sec:Acceleration_algorithm} introduces the proposed multi-cut decomposition algorithm to solve the action selection procedure. An application about multi-facility capacity investment problems, along with the problem analysis and a comprehensive numerical studies, are presented in Section \ref{sec:MCIP}. In the last section, we conclude the major findings and arguments of this study. All of our proofs may be found in the \hyperref[sec:appendix]{Appendix}.

	\section{Related Literature}
	
	\label{sec:Related_Work} 
	
	Sequential decision problems with finite and vector-valued action spaces can be found in capacity expansion problems of engineering systems, where the capacity is usually quantified by the number of production lines, machines, and facilities, which is discrete \citep{huang_value_2009, zhao_decision_2018, sun_stochastic_2015}. Similar applications can also be found in multi-item inventory control models \citep{vairaktarakis_robust_2000, topan_exact_2010}.  
	
	Sequential decision problems can be modeled as MDP and solved via ADP if the model has high-dimensional state or action space. The basic idea of ADP is to approximate the value functions or $Q$ functions via state aggregation \citep{bertsekas_dynamic_2012}, basic functions approximation \citep{munos_finite-time_2008,antos_fitted_2008,powell_approximate_2011}, or sample average approximation \citep{haskell_empirical_2016}. Many of the ADPs in the literature assumed that the action space of the MDP is finite, so that we can enumerate all possible actions in the action selection procedure. However, it can be extensively time-consuming to find the optimal action via enumeration when the action space is high-dimensional. Policy gradient algorithms have been proposed to find the optimal action for problems with continuous and high-dimensional action spaces \citep{lillicrap_continuous_2015,schulman_high-dimensional_2015}, but this type of methods may not work for derivative-free problems such as MCIPs with discrete capacity. Other methods to handle with the derivative-free problems include direct policy search\textemdash the policy of the problem is approximated, for example, via NNs\textemdash and the optimization of such policy is done by a black-box optimization \citep{hu_sequential_2017}. \citet{jain_simulation-based_2010} approximated the policies via sampling and proposed a simulation-based optimization framework for policy improvements. However, the required number of samples of this type of methods could be large. 
	
	We note that multi-stage stochastic programming is another method that can solve sequential decision problems with finite and vector-valued decisions. This method usually models the evolution of the uncertain parameters as a scenario tree, but this may have limitations as the size of the tree grows exponentially with the number of stages \citep{shapiro_lectures_2009}. To address this issue, stochastic dual dynamic integer programming is proposed, and it combines the sampling technique with a nested decomposition scheme to find the optimal policy \citep{zou_stochastic_2018}. Another stream of methods is the decision rule-based methods, where the policy of the sequential decision problems are approximated via a specific class of functions, including linear \citep{kuhn_primal_2011}, piecewise linear \citep{georghiou_generalized_2015}, binary functions \citep{bertsimas_binary_2018}, or nonparametric approach based on Gaussian processes \citep{defourny_scenario_2013}. However, these methods usually require the objective functions of the models being convex or linear. \cite{zhao_decision_2018} proposed an if-then decision rule to solve a capacity expansion problem with a non-convex objective function and vector-valued decisions, but the optimality of the solution is not guaranteed.

	\section{\label{sec:Preliminaries}Preliminaries}
	
	We consider discounted finite horizon Markov decision processes (MDP) with the following ingredients:
	\begin{itemize}
		\item planning horizon $\mathbb{T} \triangleq \{1, 2, \ldots, T\}$ with $T < \infty$ (with $\mathbb{T}\backslash\left\{ T\right\} \triangleq \{1, 2, \ldots, T-1\}$);
		\item state space $\mathbb{X}\subset\mathbb{R}^{N_{1}}$ (possibly continuous);
		\item action space $\mathbb{A}\subset\mathbb{R}^{N_{2}}$ (vector-valued and finite, but possibly very large);
		\item exogenous uncertainty $\Xi\subset\mathbb{R}^{N_{3}}$ (possibly continuous);
		\item state transition function $f : \mathbb{X} \times \mathbb{A} \times \Xi \rightarrow \mathbb{X}$;
		\item reward functions $\{r_t\}_{t \in \mathbb{T}} : \mathbb{X} \times \mathbb A \rightarrow \mathbb R$;
		\item discount factor $\gamma\in\left(0,1\right)$.
	\end{itemize}
	We assume the initial state $x_1$ is known. At time $t \in \mathbb{T}$, the current state is $x_{t}\in\mathbb{X}$, the current action is $a_{t}\in\mathbb{A}$, and the reward is $r_{t}\left(x_{t},a_{t}\right)$. Let $\xi_{t}\in\Xi$ denote the exogenous information that arrives at the end of period $t$ (which is independent of $x_{t}$ and $a_{t}$) and drives the next state transition to state $x_{t+1}=f\left(x_{t},a_{t},\xi_{t}\right)$.
	
	We gather our main assumptions below for easy reference.
	
	\begin{assumption}
		\label{assu:MDP_finite_and_discounted}The state space $\mathbb{X}$ and the support $\Xi$ are closed and bounded. The action space satisfies $\mathbb{A} \subset \times_{n=1}^{N_2} \left[ 0, \bar{a}_n \right]$. The reward functions satisfy $\left|r_{t}\left(x,a\right)\right|\le r_{\max}<\infty$ for all $(x,\,a) \in \mathbb{X} \times \mathbb A$ and $t \in \mathbb{T}$.
	\end{assumption}
	\noindent
	Our assumption that the actions are non-negative vectors is without loss of generality, since we can suitably translate any finite set of actions. We further assume that the reward function is Lipchitz in $x_t\in\mathbb{X}$.
	
	\begin{assumption}
		\label{assu:MDP_lip_reward}There exists $L_{r}<\infty$ such that $\left|r_t\left(x^{1}, a\right)-r_t\left(x^{2}, a\right)\right|\le L_{r}\left\Vert x^{1}-x^{2}\right\Vert $
		for all $x^{1},x^{2}\in\mathbb{X}$, $a\in\mathbb{A}$, and $t\in\mathbb{T}$.
	\end{assumption}
	
	Next we assume that the state transition function is affine in $a_{t}$, which is the case for many problems (e.g. inventory control, portfolio optimization, and capacity investment problems). This special structure plays a material role in our upcoming algorithm.
	
	\begin{assumption}
		\label{assu:MDP_linear_transition}The state transition function $f\left(\cdot\right)$
		is affine in $a_{t}$, i.e., $f\left(x_{t},a_{t},\xi_{t}\right) =  A\left(x_{t},\xi_{t}\right)+B\left(x_{t},\xi_{t}\right)a_{t}$ where $A:\mathbb{X}\times\Xi\mapsto\mathbb{R}^{N_{1}}$ and $B:\mathbb{X}\times\Xi\mapsto\mathbb{R}^{N_{1}\times N_{2}}$. 
	\end{assumption}
	
	The following assumption is for simplicity; we can do without it at the expense of more complicated notation.
	
	\begin{assumption}
		\label{assu:exogenous}The exogenous information $\{\xi_t\}_{t=1}^{T-1}$ is i.i.d. 
	\end{assumption}
	\noindent
	Under Assumption \ref{assu:exogenous}, there is a transition kernel $P(\cdot \, \vert \, x,\,a)$ defined by
	$$P(U \, \vert \, x,\,a) \triangleq \text{Pr}\left( \{f(x,\,a,\,\xi) \in U\} \right),
	$$
	for all measurable subsets $U \subset \mathbb{X}$ (this transition kernel is stationary because $\{\xi_t\}_{t=1}^{T-1}$ are i.i.d.).
	
	\begin{assumption}
		\label{assu:MDP_lip_transition}The transition kernel
		$P\left(\cdot\left|x,a\right.\right)$ admits a uniformly bounded
		density, i.e., there exists $L_{P}<\infty$ such that $\left|P\left(x^{1}\left|x,a\right.\right)-P\left(x^{2}\left|x,a\right.\right)\right|\le L_{P}\left\Vert x^{1}-x^{2}\right\Vert $
		for all $x^{1},x^{2}\in\mathbb{X}$ and $(x,a)\in \mathbb{X} \times \mathbb{A}$.
	\end{assumption}
	
	Let $\Pi$ be the set of all deterministic Markov policies $\pi = (\pi_t)_{t \in \mathbb{T}}$ where $\pi_t:\mathbb{X}\mapsto\mathbb{A}$ for all $t\in \mathbb{T}$. We seek a policy $\pi\in\Pi$ that maximizes the cumulative reward:
	\begin{equation}\label{MDP}
	    \sup_{\pi\in\Pi}\mathbb{E}^\pi \left\{ \sum_{t=1}^{T}\gamma^{t-1} r_{t}\left(x_{t},a_{t}\right)\right\} .
	\end{equation}
	We denote an optimal policy as $\pi^{*} = (\pi_t^{*})_{t \in \mathbb{T}}$ (it is not necessarily unique).
	
	Now we characterize the DP equations for Problem \eqref{MDP}. Let $V_{t}:\mathbb{X}\mapsto\mathbb{R}$
	denote the time $t \in \mathbb{T}$ value function. We recall the initial state $x_1$ is known, so the DP equations are:
	\begin{subequations}\label{eq:ValueFunction}
		\begin{align}
    		V_{T}\left(x\right)=\, & \max_{a\in\mathbb{A}} \left\{r_{T}\left(x,a\right)\right\},\,\forall x \in \mathbb{X},\\
    		V_{t}\left(x\right)=\, & \max_{a\in\mathbb{A}} \left\{r_{t}\left(x,a\right)+\gamma\mathbb{E}\left[V_{t+1}\left(f(x, a, \xi_{t})\right)\right]\right\},\,\forall x \in \mathbb{X},\, \forall t \in \mathbb{T}\backslash\left\{ 1,\, T\right\},\\
    		V_{1}\left(x_1\right)=\, & \max_{a\in\mathbb{A}} \left\{r_{1}\left(x_1,a\right)+\gamma\mathbb{E}\left[V_{2}\left(f(x_1, a, \xi_{1})\right)\right]\right\}.
		\end{align}
	\end{subequations}
	Theoretically speaking, we can solve Eq. \eqref{eq:ValueFunction}
	by value iteration.

	\section{Neural Network-based Fitted Value Iteration Algorithm} \label{sec:NNFVI}
	
	When the state space $\mathbb{X}$ is large or continuous, it is impossible to exactly solve Eq. \eqref{eq:ValueFunction} because of the curse of dimensionality (see \cite{powell_approximate_2011}). To respond to this challenge, we will use FVI (see \cite{munos_finite-time_2008}), which is a scheme for approximately solving Eq. \eqref{eq:ValueFunction} based on randomly sampling from the state space, estimating the value function at the sampled states, and then doing function fitting to generate an approximate value function on the entire state space.
	
	FVI depends on two main ingredients: the value function approximation architecture and the state space sampling distribution. Let $\mathbb{Y}\subset\mathbb{R}^{N_4}$ be a set of adjustable parameters which control the value function approximation. For each $y \in \mathbb{Y}$, let $\tilde{V}_{t}\left(\cdot;y\right)$ for all $t \in \mathbb{T}$ be parametrized approximate value functions on $\mathbb{X}$ (e.g. linear basis functions, polynomial basis functions, radial basis functions, etc.). Next we introduce the state space sampling distribution $\mu$, a probability distribution on the state space $\mathbb{X}$. In each iteration, we will randomly generate an i.i.d. sample of size $S_1$ indexed by $\mathcal{S}_{1}\triangleq\left\{ 1,\ldots,S_{1}\right\}$.
	
	Let us describe FVI starting with the terminal period $t = T$. First, a sample $\left(x_{T}^{s}\right)_{s\in\mathcal{S}_{1}}$ is drawn from $\mathbb{X}$ according to $\mu$, and the values $\hat{V}_{T}\left(x_{T}^{s}\right)$ are computed exactly according to Eq. \eqref{eq:ValueFunction} for $t = T$ (they can only be computed exactly in the terminal period when there is no expected future value). Then, our data set consists of the sampled states $\left(x_{T}^{s}\right)_{s\in\mathcal{S}_{1}}$ and the Bellman estimates $\left(x_{T}^{s},\hat{V}_{T}\left(x_{T}^{s}\right)\right)_{s\in\mathcal{S}_{1}}$.
	
	We want to find the best fit to the data $\left(x_{T}^{s},\hat{V}_{T}\left(x_{T}^{s}\right)\right)_{s\in\mathcal{S}_{1}}$ within some tractable class of functions. In period $T$, $\tilde{V}_{T}\left(\cdot;y\right)$ can be trained by solving the regression problem: 
	\begin{equation}\label{eq:FVI_regression_terminal}
	    \hat{y}_{T} \in \arg\min_{y \in \mathbb{Y}}\frac{1}{S_{1}}\sum_{s\in\mathcal{S}_{1}}\left(\tilde{V}_{T}\left(x_{T}^{s};y\right)-\hat{V}_{T}\left(x_{T}^{s}\right)\right)^{2}+\frac{\beta}{2}y^{\top}y,
	\end{equation}
	where $\beta>0$ is a regularization parameter. Regularizing the objective function can sometimes improve the generalization of the function fitting and avoid overfitting.
	
	After solving Problem \eqref{eq:FVI_regression_terminal} for period $T$, the algorithm proceeds inductively. In period $t\in\mathbb{T}\backslash\left\{ T\right\} $, we generate $S_{1}$ samples of the state space from $\mu$, denoted by $\left(x_{t}^{s}\right)_{s\in\mathcal{S}_{1}}$. Then, we need to simulate state transitions for each of these states, for every possible action $a_t \in \mathbb{A}$. For each state sample $x_{t}^{s}$, we generate $\xi_{t}^{s^{\prime}}$ via Monte Carlo simulation, and then for any $a_t \in \mathbb{A}$ we get a sample $x_{t+1}^{ss^{\prime}}\left(a_{t}\right)=f\left(x_{t}^{s},a_{t},\xi_{t}^{s^{\prime}}\right)$ of the next state visited. We let $\mathcal{S}_{2}\triangleq\left\{ 1,\ldots,S_{2}\right\}$ index these sampled state transitions such that $s^{\prime}\in\mathcal{S}_2$. Then, we can calculate $\hat{V}_{t}\left(x_{t}^{s}\right)$ using the trained function $\tilde{V}_{t+1}\left(\cdot;\hat{y}_{t+1}\right)$, and compute $\tilde{V}_{t}\left(\cdot;\hat{y}_{t}\right)$ by solving Problem \eqref{eq:FVI_regression}. The details of the algorithm are summarized in Algorithm \ref{alg:The-FVI-algorithm}. 
	
	\begin{algorithm}
		\caption{\label{alg:The-FVI-algorithm}The FVI algorithm}
		
		\begin{itemize}
			\item Step 0: Initialize $\mathcal{S}_{1}$, $\mathcal{S}_{2}$, initial state $x_{1}$, and the MDP parameters. Set $t\leftarrow T$.
			\item Step 1: Draw samples $\left(x_{t}^{s}\right)_{s\in\mathcal{S}_{1}}$ independently from the state space $\mathbb{X}$ and generate samples $\xi_{t}^{s^{\prime}}$ to calculate future transitions $x_{t+1}^{ss^{\prime}}\left(a_{t}\right)=f\left(x_{t}^{s},a_{t},\xi_{t}^{s^{\prime}}\right)$ for all $s\in\mathcal{S}_{1},s^{\prime}\in\mathcal{S}_{2}$.
			\item Step 2: Compute for $s\in\mathcal{S}_{1}$,
			\begin{equation}
			\hat{V}_{t}\left(x_{t}^{s}\right)=\begin{cases}
			\max\limits _{a_{t}\in\mathbb{A}}\left[r_{t}\left(x_{t}^{s},a_{t}\right)+\frac{\gamma}{S_{2}}\sum\limits _{s^{\prime}\in\mathcal{S}_{2}}\tilde{V}_{t+1}\left(x_{t+1}^{ss^{\prime}}\left(a_{t}\right);\hat{y}_{t+1}\right)\right], & 1< t<T,\\
			\max\limits_{a_{t}\in\mathbb{A}}\left[r_{t}\left(x_{t}^{s},a_{t}\right)\right], & t=T.
			\end{cases}\label{eq:FVI_step2}
			\end{equation}
			
			\item Step 3: Given samples $\left(x_{t}^{s},\hat{V}_{t}\left(x_{t}^{s}\right)\right)_{s\in\mathcal{S}_{1}}$,
			fit the approximate value functions 
			\begin{equation}
			\hat{y}_{t}\in\arg\min_{y \in \mathbb{Y}}\frac{1}{S_{1}}\sum_{s\in\mathcal{S}_{1}}\left(\tilde{V}_{t}\left(x_{t}^{s};y\right)-\hat{V}_{t}\left(x_{t}^{s}\right)\right)^{2}+\frac{\beta}{2}y^{\top}y.\label{eq:FVI_regression}
			\end{equation}
			
			\item Step 4: If $t>1$, set $t\leftarrow t-1$ and go to Step	1; otherwise, generate samples $\xi_2^{s^\prime}$ to calculate future transitions $x_2^{s^{\prime}}(a_1)=f\left(x_{1},a_{1},\xi_{2}^{s^{\prime}}\right)$, return
			\[
			\hat{V}_{1}\left(x_{1}\right)=\max_{a_{1}\in\mathbb{A}}\left[r_{1}\left(x_{1},a_{1}\right)+\frac{\gamma}{S_{2}}\sum_{s^{\prime}\in\mathcal{S}_{2}}\tilde{V}_{2}\left(x_{2}^{s^{\prime}}\left(a_{1}\right);\hat{y}_{2}\right)\right],
			\]
			and terminate the algorithm.
		\end{itemize}
	\end{algorithm}
	
	There are two questions that need to be answered about the above FVI algorithm. First, what approximation architecture $\left\lbrace \tilde{V}_{t}\left(\cdot;y\right) \right\rbrace_{t \in \mathbb{T}}$ should we use for the value functions? Second, the action selection problem in Eq. \eqref{eq:FVI_step2} optimizes over a multi-dimensional action space $\mathbb{A}$. It is difficult to solve this problem via direct enumeration, so how can we do this optimization? We will answer both questions. Next, in Subsection \ref{subsec:NNFVI} we choose a value function approximation architecture based on NNs. Then, in Section \ref{sec:Acceleration_algorithm} we introduce a decomposition algorithm to accelerate the action selection in Problem \eqref{eq:FVI_step2}.
	
	\subsection{\label{subsec:NNFVI}Two-layer Neural Network with ReLU}
	
	There is a well-known exploration vs. exploitation trade-off in general function approximation problems. One can use a rich functional family, such as deep NNs, to approximate the value functions in Eq. \eqref{eq:ValueFunction}. However, the training and inference may then be time-consuming given that we need to solve the action selection problem $\left(T-1\right) \times S_1+1$ times. In this paper, we propose an FVI, where the value functions are approximated via a two-layer NN with rectified linear units (ReLU). As we will see, this approximation is powerful enough to reach any desired accuracy in our problem, as long as sufficient neurons are provided. We will see also that this approximation leads to efficient optimization in the action selection problem.
	
\begin{figure}
		\centering
		{\includegraphics[scale=0.3]{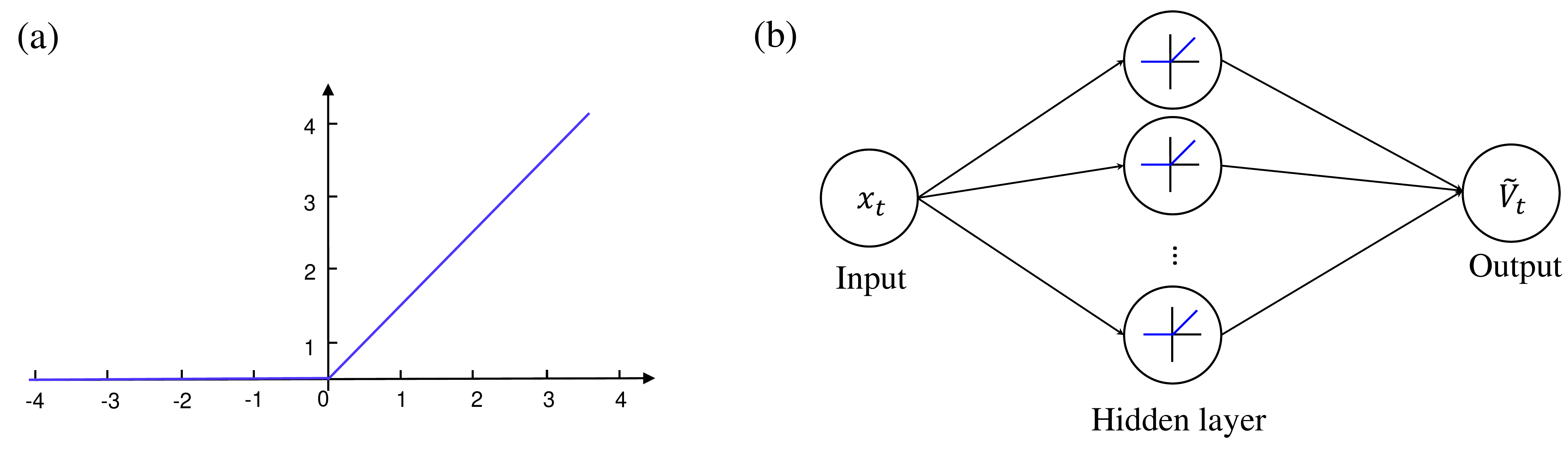}}
		\caption{\label{fig:Neural Network}(a) the plot of ReLU and (b) the architecture of the two-layer NN}
		{}
	\end{figure}

	A two-layer NN consists of inputs (the state $x_{t}\in\mathbb{X}$), one hidden layer for
	intermediate computations, and a single output (the approximate value $\tilde{V}_{t}\left(\cdot\right)$ of the value function). Let $\mathcal{J}\triangleq\left\{ 1,\ldots,J\right\} $
	index the neurons in the hidden layer of our network, so the two-layer NN is
	\[
	\Gamma\left(x\right)=\sum_{j\in\mathcal{J}}w_{j}\Psi_{j}\left(u_{j}^{\top}x+u_{0j}\right)+w_{0j},
	\]
	where $\left(u_{j},u_{0j}\right)$ and $\left(w_{j},w_{0j}\right)$
	are the adjustable weights of the input layer and the hidden layer
	respectively, and $\Psi_{j}\left(\cdot\right)$ is the activation
	function for neuron $j \in \mathcal{J}$. We choose the activation functions $\Psi_{j}\left(\cdot\right)$ to
	be ReLU, i.e.,
	\[
	\Psi_{j}\left(u_{j}^{\top}x+u_{0j}\right)=\max\left\{ u_{j}^{\top}x+u_{0j},0\right\},\,\forall j \in \mathcal{J},
	\]
	which are themselves piecewise linear. The architecture of the NN is shown in Figure \ref{fig:Neural Network}. 
	
	The adjustable weights for the input of neuron $j$ in $\tilde{V}_{t}\left(\cdot\right)$
	include $u_{0jt}$ and a vector of $N_{1}$ (the dimension of $\mathbb{X}$) adjustable coefficients
	$u_{jt}=\left(u_{1jt},\ldots,u_{N_{1}jt}\right)$ for all $j\in\mathcal{J}$.
	Let $u_{t}=\left(u_{0jt},u_{jt}\right)_{j\in\mathcal{J}}$ denote the
	vector of all adjustable weights for the input layer and let
	$w_{t}=\left(w_{0t},w_{1t},\ldots,w_{Jt}\right)$
	denote the vector of the adjustable weights for the output of
	the hidden layer. We then let $y_t \triangleq \left(u_{t},w_{t}\right)$ succinctly denote all parameters for the period $t$ value function approximation. The NN at time $t\in\mathbb{T}$ is then a function of the input $x_{t}\in\mathbb{X}$
	and the adjustable weights $y_t$:
	\begin{equation}
	\tilde{V}_{t}\left(x_t;y_t\right)=\sum_{j\in\mathcal{J}}w_{jt}\max\left\{ u_{jt}^{\top}x_t+u_{0jt},0\right\} +w_{0t},\ \ \forall t\in\mathbb{T}.\label{eq:Neural_Network}
	\end{equation}
	We use the specific approximation in Eq. \eqref{eq:Neural_Network}
	in Step 3 of Algorithm \ref{alg:The-FVI-algorithm}. From now on, we refer to our algorithm as NN-based FVI (NN-FVI). 
	
	\subsection{Consistency of the algorithm}
	
	In this subsection, we show that our customized NN is powerful
	enough to approximate the value funtions of our MDP. Specifically, we show
	that NN-FVI is consistent.
	
	First, we show that the value functions $\left\{ V_{t}\right\} _{t\in\mathbb{T}}$
	are Lipschitz continuous under our assumptions.
	\begin{proposition}
		\label{prop:Lip_Value_Function}Suppose Assumptions \ref{assu:MDP_finite_and_discounted}--\ref{assu:MDP_lip_transition} hold, then
		$V_{t}\left(\cdot\right)$ is $L_{v}-$Lipschitz for all $t\in\mathbb{T}$ (i.e., there exists
		$L_{v}<\infty$ such that $\left|V_{t}\left(x\right)-V_{t}\left(x^{\prime}\right)\right|\le L_{v}\left\Vert x-x^{\prime}\right\Vert $
		for all $x,x^{\prime}\in\mathbb{X}$ and $t\in\mathbb{T}$).
	\end{proposition}
	\noindent
	We then define the class
	\[
	\text{Lip}\left(L_{v}\right)\triangleq\left\{ v:\mathbb{X}\mapsto\mathbb{R}\left|\left|v\left(x\right)-v\left(x^{\prime}\right)\right|\le L_{v}\left\Vert x-x^{\prime}\right\Vert ,\forall x,x^{\prime}\in\mathbb{X}\right.\right\} 
	\]
	of $L_{v}$-Lipschitz functions. The approximation
	power of a function class $\mathcal{F}$ on $\mathbb{X}$ with respect to $\text{Lip}\left(L_{v}\right)$ is measured by the inherent Bellman error (see \cite{munos_finite-time_2008}):
	\[
	\mathcal{E}_{p}\left(\text{Lip}\left(L_{v}\right),\mathcal{F}\right)\triangleq\sup_{g\in\text{Lip}\left(L_{v}\right)}\inf_{f\in\mathcal{F}}\left\Vert g-f\right\Vert _{p},
	\]
	for $1 \leq p \leq \infty$. Since $V_{t}\left(\cdot\right)\in\text{Lip}\left(L_{v}\right)$, if $\mathcal{E}_{p}\left(\text{Lip}\left(L_{v}\right),\mathcal{F}\right)$
	is small then the function class $\mathcal{F}$ contains a good approximation of
	$V_{t}\left(\cdot\right)$. If a sequence $\left\{ \mathcal{F}_{i}\right\}_{i \geq 1} $ satisfies
	$\lim_{i\rightarrow\infty}\mathcal{E}_{p}\left(\text{Lip}\left(L_{v}\right),\mathcal{F}_{i}\right)=0$,
	then it is \emph{universal }\citep{antos_value-iteration_2007}\emph{.}
	This condition implies that the inherent approximation error between
	the two function spaces converges to zero as $\left\{ \mathcal{F}_{i}\right\}_{i \geq 1} $
	becomes ``richer''.
	
	The richness of an NN is quantified by the number of layers
	and weights \citep{bartlett_nearly-tight_2019}. The approximation
	power of our customized two-layer NN increases with the number
	of neurons, and we have the following result about its approximation quality.
	\begin{lemma}
		\label{lem:Universal_Appro}The two-layer NN with ReLU activation functions
		is a universal approximator. 
	\end{lemma}
	
	Given Assumptions
	\ref{assu:MDP_finite_and_discounted} and \ref{assu:MDP_lip_transition},
	the consistency of NN-FVI follows from \citep[Corollary 4]{munos_finite-time_2008}.
	We state this result below and present its detailed proof in Appendix
	\ref{sec:Proof-of-Consistency}.
	\begin{theorem}
		\label{prop:Consistency}Suppose Assumptions \ref{assu:MDP_finite_and_discounted} and		\ref{assu:MDP_lip_transition} hold. For any $\epsilon>0$ and $\delta\in\left(0,1\right)$, there exists an integer $J_{0}$ such that for any $J\ge J_{0}$ there are $S_{1}$, $S_{2}$ that are polynomial in $1/\epsilon$, $\log(1/\delta)$, $1/(1-\gamma)$, $r_{\max}$, $\log(|\mathbb{A}|)$, and the pseudo-dimension of the NN (see \cite{haussler1995sphere,munos_finite-time_2008}) such that $\left\Vert \tilde{V}_{t}-V_{t}\right\Vert _{\infty}\le\epsilon$ holds for all $t\in\mathbb{T}$ with probability at least $1-\delta$.
	\end{theorem}
	
	\begin{remark}
		Theorem \ref{prop:Consistency} confirms the consistency of NN-FVI.
		However, the NN training problem Eq. \eqref{eq:FVI_regression}
		in Algorithm \ref{alg:The-FVI-algorithm} is non-convex. Fortunately,
		it has been demonstrated numerically that the chance of getting stuck at a poor local minimum decreases as the size of the NN increases. In addition,
		when the NN is large enough, finding the global minimum may also be
		unnecessary as it often leads to overfitting \citep{choromanska_loss_2014}.
	\end{remark}

	\section{Accelerate the Action Selection Procedure in NN-FVI}
	
	\label{sec:Acceleration_algorithm}
	
	In NN-FVI, once the approximate value function in time $t+1$ has been fitted, one needs to solve Problem \eqref{eq:FVI_step2} (the action selection problem). This problem can in principle be solved by enumerating all possible actions since $\mathbb{A}$ is finite (just the brute-force method). However, we are concerned with the situation where the action space is very large so this method is not practical for our purposes. For example, consider a capacity expansion problem with five facilities, each of which can install $\{0,1,\ldots,\,9\}$ units of capacity (e.g. production lines). The complexity of brute-force action selection is $10^{5}$ for each instance of the action selection problem, which makes the overall FVI algorithm time-consuming. To address this issue, we formulate the action selection problem as a two-stage stochastic programming problem. This transformation allows us to achieve a speed up with a specialized decomposition algorithm.
	
	\subsection{Formulation of the Action Selection Problem}
	
	For period $t \in \mathbb{T}\backslash\left\{ T\right\}$, suppose
	we have already trained the NN in time $t+1$ and computed its adjustable weights $y_{t+1}$ in Step 3 of Algorithm
	\ref{alg:The-FVI-algorithm}. The upcoming procedure is the same for each sampled state in $\mathcal{S}_{1}$, so with some abuse of notation we suppress the dependence of the coefficients and parameters on the particular sampled state. Instead, we just let $x_t$ denote a specific state sample in $\mathcal{S}_{1}$, and we write $x_{t+1}^{s} = f\left(x_{t},a_{t},\xi_{t}^{s}\right)$ for all $s\in\mathcal{S}_{2}$ to denote samples of future transitions that are generated via Monte Carlo simulation given a specific state-action pair $(x_{t},\,a_{t})$.
	
	We define $v_{t+1}\left(x_{t},a_{t};y_{t+1}\right)$ to be the empirical recourse function
	given action $a_{t}$ and $x_{t}$ and the next state samples in $\mathcal{S}_2$, defined by:
	\begin{equation}
	v_{t+1}\left(x_{t},a_{t};y_{t+1}\right)\triangleq \frac{1}{S_2}\sum_{s\in\mathcal{S}_{2}}\sum_{j\in\mathcal{J}}\left[w_{j\left(t+1\right)}\max\left\{ u_{j\left(t+1\right)}^{\top}f\left(x_{t},a_{t},\xi_{t}^{s}\right)+u_{0j\left(t+1\right)},0\right\} \right]
	.\label{eq:recourse}
	\end{equation}
	Then, Problem \eqref{eq:FVI_step2} which computes $\hat{V}_{t}\left(x_{t}\right)$ is equivalent to:
	\begin{alignat}{1}
	\max_{a_{t}\in\mathbb{A}}\  & r_{t}\left(x_{t},a_{t}\right)+\gamma v_{t+1}\left(x_{t},a_{t};y_{t+1}\right)+\gamma w_{0(t+1)}.\label{eq:OAS_obj}
	\end{alignat}
	Note that $x_t$ is a known state sample and $w_{0(t+1)}$ is a constant, so $a_t$ is the only decision variable in Problem \eqref{eq:OAS_obj}. This problem can be viewed as a two-stage stochastic programming problem with a simple recourse function: the first-stage is to determine the action $a_{t}$, and then the recourse function returns the expectation of the future reward based on the NN trained in time $t+1$.
	
	Problem \eqref{eq:OAS_obj} has two notable properties:
	\renewcommand{\labelenumi}{\arabic{enumi}.}
	\begin{enumerate}
		\item[] \textit{Property 1.} The recourse function is complete; that is, given
		any $a_{t}\in\mathbb{\mathbb{A}}$, the recourse is not empty.
		\item[] \textit{Property 2.} The recourse function can be decomposed into a sum indexed by $\mathcal{S}_2$, given any $a_t$.
	\end{enumerate}
	
	Given these properties, one may prefer to decompose Problem \eqref{eq:OAS_obj} and then update the decisions by a cut generation method (e.g. Benders decomposition), especially when $S_2$ is large. This type of algorithm is based on constructing some valid inequalities for the feasible region of Problem \eqref{eq:OAS_obj}. We recall the hypograph of the recourse function is: 
	\begin{alignat*}{1}
	    \text{hypo}_{t} & \triangleq\left\{ \left(a_{t},\eta\right)\in\mathbb{R}^{N_2+1}\left|\eta\le
	    \frac{1}{S_2}\sum_{s\in\mathcal{S}_{2}}\sum_{j\in\mathcal{J}}\left[w_{j\left(t+1\right)}\max\left\{ u_{j\left(t+1\right)}^{\top}f\left(x_{t},a_{t},\xi_{t}^{s}\right)+u_{0j\left(t+1\right)},0\right\} \right]\ \ \right.\right\}.
	\end{alignat*}
	Our procedure is based on constructing valid inequalities for the hypograph of the recourse function.
	
	\begin{definition}
		An inequality $\eta\le\phi^{\top}a_{t}+\phi_{0}$ is \emph{valid} for the set $\text{hypo}_{t}$ if 
		\[
		\text{hypo}_{t}\subset\left\{ \left(a_{t},\eta\right)\in\mathbb{R}^{N_2+1}\left\vert\eta\le\phi^{\top}a_{t}+\phi_{0}\right.\right\} .
		\]
	\end{definition}
	
	A cut generation algorithm approximates the hypograph of Problem \eqref{eq:OAS_obj} via some valid cuts iteratively. Let $m$ index the iterations of the algorithm, and let $a_{t}^{m}$ denote the optimal solution of the first-stage problem in the $\left(m-1\right)$th iteration. The approximate first-stage problem in the $m$th iteration is of the general form:
	\begin{subequations}
		\label{eq:OAS_obj_hypo}
		\begin{alignat}{2}
		\max_{a_{t}\in\mathbb{A},\,\eta\in \mathbb{R}} \  & r_{t}\left(x_{t},a_{t}\right)+\gamma \eta+\gamma w_{0\left(t+1\right)} &\\
		\text{s.t. } & \eta\le\left(\phi^{m^{\prime}}\right)^{\top}a_{t}+\phi^{m^{\prime}}_{0}, \ \ & \forall m^{\prime}=1,\ldots,m,
		\end{alignat}
	\end{subequations}
	where the inequality constraints are the cuts constructed so far. The decomposition procedure is as follows: 
	\begin{enumerate}
		\item First, we solve the first-stage problem and compute $a_{t}^{m}$
		in the $\left(m-1\right)$th iteration. 
		\item \label{enu:Step_2}Then, we construct a valid inequality given $a_{t}^{m}$,
		add it into the first-stage problem, and update the first-stage decision
		in the $m$th iteration. 
		\item The algorithm is run iteratively and the approximation of the recourse function is improved via valid inequalities from above. 
	\end{enumerate}
	
	If the recourse function is concave in $a_{t}$, then we can compute	a supporting hyperplane of the hypograph of the recourse function at any	$a_{t}\in\mathbb{A}$ and this hyperplane would be a valid inequality.	In this situation, the cut generation algorithm behaves exactly like Benders decomposition.
	
\begin{figure}
		\centering
		{\includegraphics[scale=0.3]{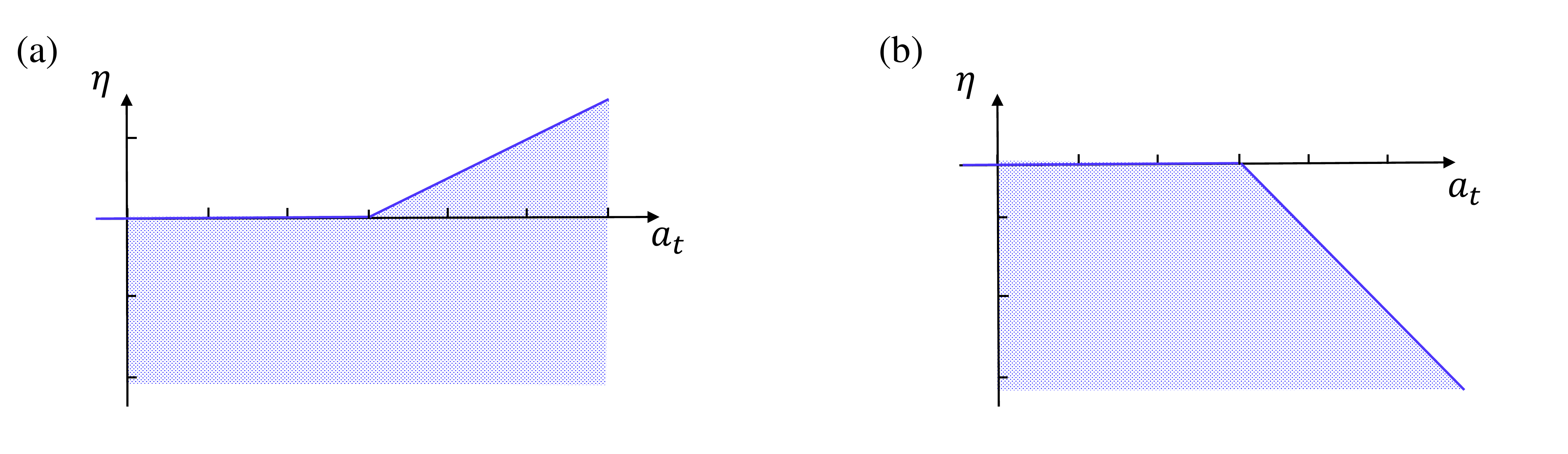}}
		\caption{\label{fig:hypograph_neurons}Hypographs of the output of neuron $j$ when (a) $w_{j\left(t+1\right)}>0$ and (b) $w_{j\left(t+1\right)}\le0$}
		{}
	\end{figure}

	However, since the activation function of our NN is ReLU, the recourse function may not be concave in $a_t\in\mathbb{A}$.	As can be seen in Fig. \ref{fig:hypograph_neurons}, the hypograph of neuron $j$ is a convex set if $w_{j\left(t+1\right)}\le0$. On the other hand, if $w_{j\left(t+1\right)}>0$, then the hypograph is not convex in which case the objective of Problem \eqref{eq:OAS_obj} may fail to be concave with respect to $a_t$. The sign of $u_{j\left(t+1\right)}$ does not effect convexity of the hypograph. We have the following property.
	
	\begin{enumerate}
		\item[] \textit{Property 3.} The recourse function in Eq. \eqref{eq:recourse} is not necessarily concave in $a_{t}\in\mathbb{\mathbb{A}}$ if there exists any neuron $j\in\mathcal{J}$ such that $w_{j\left(t+1\right)}>0$. 
	\end{enumerate}
	
	Based on these observations, we separate the neurons with positive
	and non-positive weights into different sets:
	\[
	\mathcal{J}_{t+1}^{+}\triangleq\left\{ \left.j\in\mathcal{J}\right|w_{j\left(t+1\right)}>0\right\} \ \ \text{and}\ \ \mathcal{J}_{t+1}^{-}\triangleq\left\{ \left.j\in\mathcal{J}\right|w_{j\left(t+1\right)}\le0\right\} .
	\]
	Then, we may define the hypographs
	\begin{alignat*}{1}
	\text{hypo}_{t}^{+} & \triangleq\left\{ \left(a_{t},\eta^{+}\right)\in\mathbb{R}^{N_2+1}\left|\eta^{+}\le
	\frac{1}{S_2}\sum_{s\in\mathcal{S}_{2}}\sum_{j\in\mathcal{J}_{t+1}^{+}}\left[w_{j\left(t+1\right)}\max\left\{ u_{j\left(t+1\right)}^{\top}f\left(x_{t},a_{t},\xi_{t}^{s}\right)+u_{0j\left(t+1\right)},0\right\} \right]
	\right.\right\} ,\\
	\text{hypo}_{t}^{-} & \triangleq\left\{ \left(a_{t},\eta^{-}\right)\in\mathbb{R}^{N_2+1}\left|\eta^{-}\le
	\frac{1}{S_2}\sum_{s\in\mathcal{S}_{2}}\sum_{j\in\mathcal{J}_{t+1}^{-}}\left[w_{j\left(t+1\right)}\max\left\{ u_{j\left(t+1\right)}^{\top}f\left(x_{t},a_{t},\xi_{t}^{s}\right)+u_{0j\left(t+1\right)},0\right\} \right]\right.\right\} ,
	\end{alignat*}
	corresponding to the output of the neurons with positive/non-positive
	weights, respectively. 
	
	In the following subsections, we design valid cuts for $\text{hypo}_{t}^{+}$ and $\text{hypo}_{t}^{-}$ separately, and find the optimal action via a decomposition algorithm. 
	
	\subsection{Valid Cuts for the Action Selection Problem}
	
	\subsubsection{Valid cuts for $\text{hypo}_t^-$}
	
	In this subsection, we construct valid cuts for neurons with negative weights (i.e., for $j \in \mathcal{J}_{t+1}^{-}$)
	by using gradient information. Suppose we have solved
	the first-stage problem in the $\left(m-1\right)$th iteration and
	computed the first-stage decision $a_{t}^{m}$. In Problem \eqref{eq:OAS_obj},
	the gradient of the output of neuron $j\in\mathcal{J}_{t+1}^{-}$ at point $a_{t}^{m}$
	can be computed directly, and the corresponding supporting hyperplane is a valid cut for $\text{hypo}_{t}^{-}$ since it is convex.
	
	Let $\mathbb{I}\left(\cdot\right)$
	be the indicator function for the positive reals where $\mathbb{I}\left(z\right)=1$
	if $z>0$ and $\mathbb{I}\left(z\right)=0$ otherwise. We propose the following cuts for neurons with negative
	weights:
	\begin{alignat}{1}
	\eta^{-} & \le
	\frac{1}{S_2}\sum_{s\in\mathcal{S}_{2}}\sum_{j\in\mathcal{J}_{t+1}^{-}}\left(\left(\phi_{j}^{ms}\right)^{\top}a_{t}+\phi_{0j}^{ms}\right),\ \ \forall a_{t}\in\mathbb{A},\label{eq:validcut_neg}
	\end{alignat}
	where:
	\begin{align*}
	\phi_{j}^{ms}\triangleq\, & w_{j\left(t+1\right)}u_{j\left(t+1\right)}^{\top}B\left(x_{t},\xi_{t}^{s}\right)\mathbb{I}\left(u_{j\left(t+1\right)}^{\top}f\left(x_{t},a_{t}^{m},\xi_{t}^{s}\right)+u_{0j\left(t+1\right)}\right),\,\forall j\in\mathcal{J}_{t+1}^{-},\,s\in\mathcal{S}_{2},\\
	\phi_{0j}^{ms}\triangleq\, & w_{j\left(t+1\right)}\left(u_{j\left(t+1\right)}^{\top}A\left(x_{t},\xi_{t}^{s}\right)+u_{0j\left(t+1\right)}\right)\mathbb{I}\left(u_{j\left(t+1\right)}^{\top}f\left(x_{t},a_{t}^{m},\xi_{t}^{s}\right)+u_{0j\left(t+1\right)}\right), \,\forall j\in\mathcal{J}_{t+1}^{-},\,s\in\mathcal{S}_{2}.
	\end{align*}
	Note that $\phi_{j}^{ms}\in\mathbb{R}^{N_{2}}$ and $\phi_{0j}^{ms}\in\mathbb{R}$
	depend on the index $s\in\mathcal{S}_{2}$. 
	
	\begin{proposition} 
		\label{prop:valid_inquality_neg} Cut \eqref{eq:validcut_neg} is a valid inequality for $\text{hypo}_{t}^{-}$.
	\end{proposition}
	
	\subsubsection{Valid cuts for $\text{hypo}_t^+$}
	
	To derive a valid cut for neurons $j\in\mathcal{J}_{t+1}^{+}$, we
	need to find a supporting hyperplane for $\text{hypo}_{t}^{+}$, which is a non-convex set. We propose the following inequality
	\begin{equation}
	\eta^{+}\le\frac{1}{S_{2}}\sum_{s\in\mathcal{S}_{2}}\sum_{j\in\mathcal{J}_{t+1}^{+}}\left(\sum_{n=1}^{N_{2}}\theta_{nj}^{s}a_{nt}+\theta_{0j}^{s}\right),\ \ \forall a_{t}\in\mathbb{A},\label{eq:validcut_pos}
	\end{equation}
	which is constructed in the following way. First, we write
	\[
	u_{j\left(t+1\right)}^{\top}f\left(x_{t},a_{t},\xi_{t}^{s}\right)+u_{0j\left(t+1\right)}=\sum_{n=1}^{N_2}\Gamma_{njs}^{\left(1\right)}a_{nt}+\Gamma_{js}^{\left(2\right)},\ \ \forall s\in\mathcal{S}_{2},
	\]
	where $\Gamma_{njs}^{\left(1\right)}\triangleq \left[u_{j\left(t+1\right)}^{\top}B\left(x_{t},\xi_{t}^{s}\right)\right]_n$ and $\Gamma_{js}^{\left(2\right)}\triangleq u_{j\left(t+1\right)}^{\top} A\left(x_{t},\xi_{t}^{s}\right)+u_{0j\left(t+1\right)}$, and we define sets
	\[
	\mathcal{\mathbf{N}}_{js}^{+}\triangleq \left\{  n=1,\ldots,N_{2}\left|\Gamma_{njs}^{\left(1\right)} \ge0 \right.\right\}
	\text{ and }
	\mathbf{N}_{js}^{-}\triangleq\left\{ n=1,\ldots,N_{2}\left|\Gamma_{njs}^{\left(1\right)}<0\right.\right\}. 
	\]
	Because the function $u_{j\left(t+1\right)}^{\top}f\left(x_{t},a_{t},\xi_{t}^{s}\right)+u_{0j\left(t+1\right)}$ is linear in $a_t$, its minimum/maximum values can be computed explicitly and are $\sum_{n\in\mathbf{N}_{js}^{-}}\Gamma_{njs}^{\left(1\right)}\bar{a}_{n}+\Gamma_{js}^{\left(2\right)}$ and $\sum_{n\in\mathbf{N}_{js}^{+}}\Gamma_{njs}^{\left(1\right)}\bar{a}_{n}+\Gamma_{js}^{\left(2\right)}$, respectively. We denote
	
	\[
	\vartheta_{j}^{s}\triangleq\frac{\sum_{n\in\mathbf{N}_{js}^{+}}\Gamma_{njs}^{\left(1\right)}\bar{a}_{n}+\Gamma_{js}^{\left(2\right)}}{\sum_{n=1}^{N_{2}}\left|\Gamma_{njs}^{\left(1\right)}\right|\bar{a}_{n}}.
	\]
	We then define for all $n=1,\ldots,N_{2}$ the constants
	\[
	\theta_{nj}^{s}\triangleq\left\{ \begin{array}{cc}
	w_{j\left(t+1\right)}\Gamma_{njs}^{\left(1\right)}, & \text{if }\sum_{n\in\mathbf{N}_{js}^{-}}\Gamma_{njs}^{\left(1\right)}\bar{a}_{n}+\Gamma_{js}^{\left(2\right)}>0,\\
	0, & \text{if }\sum_{n\in\mathbf{N}_{js}^{+}}\Gamma_{njs}^{\left(1\right)}\bar{a}_{n}+\Gamma_{js}^{\left(2\right)}<0,\\
	w_{j\left(t+1\right)}\Gamma_{njs}^{\left(1\right)}\vartheta_{j}^{s}, & \text{otherwise,}
	\end{array}\right.\ \
	\]
	and
	\[
	\theta_{0j}^{s}\triangleq\left\{ \begin{array}{cc}
	w_{j\left(t+1\right)}\Gamma_{js}^{\left(2\right)}, & \text{if }\sum_{n\in\mathbf{N}_{js}^{-}}\Gamma_{njs}^{\left(1\right)}\bar{a}_{n}+\Gamma_{js}^{\left(2\right)}>0,\\
	0 & \text{if }\sum_{n\in\mathbf{N}_{js}^{+}}\Gamma_{njs}^{\left(1\right)}\bar{a}_{n}+\Gamma_{js}^{\left(2\right)}<0,\\
	-w_{j\left(t+1\right)}\sum_{n\in\mathbf{N}_{js}^{-}}\vartheta_{j}^{s}\Gamma_{njs}^{\left(1\right)}\bar{a}_{n}, & \text{otherwise},
	\end{array}\right.
	\]
	to complete the construction of Cut \eqref{eq:validcut_pos}. This cut can be computed easily because of its linear structure. In addition, Cut \eqref{eq:validcut_pos} is independent of the current $a_t^m$ and so we only need to compute it once for a specific state sample in $\mathcal{S}_1$. 
	
	\begin{proposition}
		\label{prop:valid_inquality_pos}
		Cut \eqref{eq:validcut_pos} is a valid inequality for $\text{hypo}_{t}^{+}$.
	\end{proposition}
	
	\begin{remark}
		When training an NN, one may prefer to use leaky ReLU, instead of ReLU, as activation functions to avoid the vanishing gradient issue (or the so-called ``dying ReLU'' problem). Our valid cuts are also applicable to NNs with leaky ReLU after some minor modification.
	\end{remark}
	
	We sum up Cuts \eqref{eq:validcut_neg} (for negative neurons) and \eqref{eq:validcut_pos} (for positive neurons) to immediately get the following result.
	\begin{corollary}
	    \label{coro:valid_cut_all}
		The following cut is a valid inequality for $\text{hypo}_t$:
		\begin{equation}
		\eta\le\frac{1}{S_2}\sum_{s\in\mathcal{S}_{2}}\left[\sum_{j\in\mathcal{J}_{t+1}^{-}}\left(\sum_{n=1}^{N_2} \phi_{nj}^{ms}a_{nt}+\phi_{0j}^{ms}\right)
		+ \sum_{j\in\mathcal{J}_{t+1}^{+}}  \left( \sum_{n=1}^{N_2} \theta_{nj}^{s}a_{nt}+\theta_{0j}^{s} \right)
		\right], \ \forall a_t \in \mathbb{A}. \label{eq:validcut_all}
		\end{equation}
		
	\end{corollary}
	
	\subsubsection{Integer Optimality Cuts}
	
	In each iteration of the decomposition algorithm, we add Cuts \eqref{eq:validcut_all} to the first-stage problem. These cuts are valid, but may not be tight. 
	To ensure global optimality, we also introduce the \emph{integer optimality cuts} from \citep{laporte_integer_1993} into our algorithm. When the first-stage decision $a_{t}^{m}$ is solved in the $\left(m-1\right)$th iteration, an integer optimality cut is constructed and added to the first-stage problem. This cut exactly recovers the recourse function value corresponding to $a_{t}^{m}$ if the first-stage decision is equal to $a_{t}^{m}$, and it recovers an upper bound otherwise. 
	
	To formulate integer optimality cuts, let $\bar{\eta}$ be an upper bound on the value of the recourse function. Also define the function
	$\zeta^{m}\left(a_{t}\right)$ where $\zeta^{m}\left(a_{t}\right)=0$ if $a_{t}=a_{t}^{m}$ and $\zeta^{m}\left(a_{t}\right)\ge1$
	otherwise (we may need to transform the general integer variables $a_{t}$ into binary variables to properly formulate $\zeta^{m}\left(\cdot\right)$, the details are in the \hyperref[sec:appendix]{Appendix}). From \citep{laporte_integer_1993}, we then obtain the following integer optimality cuts:
	\[
	\eta\le v_{t+1}\left(x_{t},a_{t}^{m};y_{t+1}\right)+\zeta^{m}\left(a_{t}\right)\left[\bar{\eta}-v_{t+1}\left(x_{t},a_{t}^{m};y_{t+1}\right)\right],
	\]
	where $v_{t+1}\left(\cdot\right)$ is given by Eq. \eqref{eq:recourse}.
	This cut recovers $v_{t+1}\left(x_{t},a_{t}^{m};y_{t+1}\right)$
	if $\zeta^{m}\left(a_{t}\right)=0$, i.e., $a_{t}=a_{t}^{m}$; otherwise, it
	recovers an upper bound.
	
	
	\subsection{A Multi-Cuts Decomposition Algorithm}
	
	We now present our complete multi-cuts decomposition (MCD) algorithm for solving Problem \eqref{eq:OAS_obj}. In our algorithm, Cuts \eqref{eq:validcut_neg} and \eqref{eq:validcut_pos} exploit the structure of the recourse function, and the integer optimality cuts ensure global convergence. We continue to let $m$ index the iterations of our decomposition algorithm and $a_{t}^{m}$ denote the optimal solution of the	first-stage problem in the $\left(m-1\right)$th iteration. The first-stage problem in the $m$th iteration of the MCD algorithm is
	
	\vspace{-0.5cm}
	{\small{}
	\begin{alignat}{2}
		\tag{{\ensuremath{\mathcal{FP}}}}\max_{a_{t}\in\mathbb{A},\eta\in\mathbb{R}}\  & r_{t}\left(x_{t},a_{t}\right)+\gamma\eta +\gamma w_{0\left(t+1\right)} \label{eq:FP}\\
		\text{s.t. } & \eta\le v_{t+1}\left(x_{t},a_{t}^{m^{\prime}};y_{t+1}\right)+\zeta^{m^{\prime}}\left(a_{t}\right)\left[\bar{\eta}-v_{t+1}\left(x_{t},a_{t}^{m^{\prime}};y_{t+1}\right)\right],\ & \forall m^{\prime}=1,\ldots,m,\label{eq:LLCuts_MP}\\
		& \eta\le\frac{1}{S_2}\sum_{s\in\mathcal{S}_{2}}\left[\sum_{j\in\mathcal{J}_{t+1}^{-}}\left(\sum_{n=1}^{N_2}  \phi_{nj}^{m^{\prime}s}a_{nt}+\phi_{0j}^{m^{\prime}s}\right)
		+ \sum_{j\in\mathcal{J}_{t+1}^{+}} \left(\sum_{n=1}^{N_2} \theta_{nj}^{s}a_{nt}+\theta_{0j}^{s} \right)
		\right], & \forall m^{\prime}=1,\ldots,m. \label{eq:validcut_all2}
	\end{alignat}
	}There are two types of cuts in the above formulation. Cuts \eqref{eq:LLCuts_MP} are the integer optimality cuts generated up to iteration $m$, and Cuts \eqref{eq:validcut_all2} follow directly from Cuts \eqref{eq:validcut_all}.
	
	Our algorithm proceeds by adding Cuts \eqref{eq:LLCuts_MP}--\eqref{eq:validcut_all2} to Problem \eqref{eq:FP} simultaneously to update the first-stage decision. An upper bound can be obtained by solving Problem \eqref{eq:FP} in the $m$th iteration, and a lower bound can be obtained from the objective value of the best-found decision up to iteration $m$. One can therefore terminate the algorithm when the difference between the upper and lower bounds falls below a preset tolerance. The finite convergence of the MCD algorithm can be established as follows. 
	\begin{theorem}
		\label{prop:finite-conveg}The MCD algorithm yields a globally optimal		solution in a finite number of iterations. 
	\end{theorem}
	
\begin{figure}[!tph]
		\begin{centering}
			\includegraphics[scale=0.55]{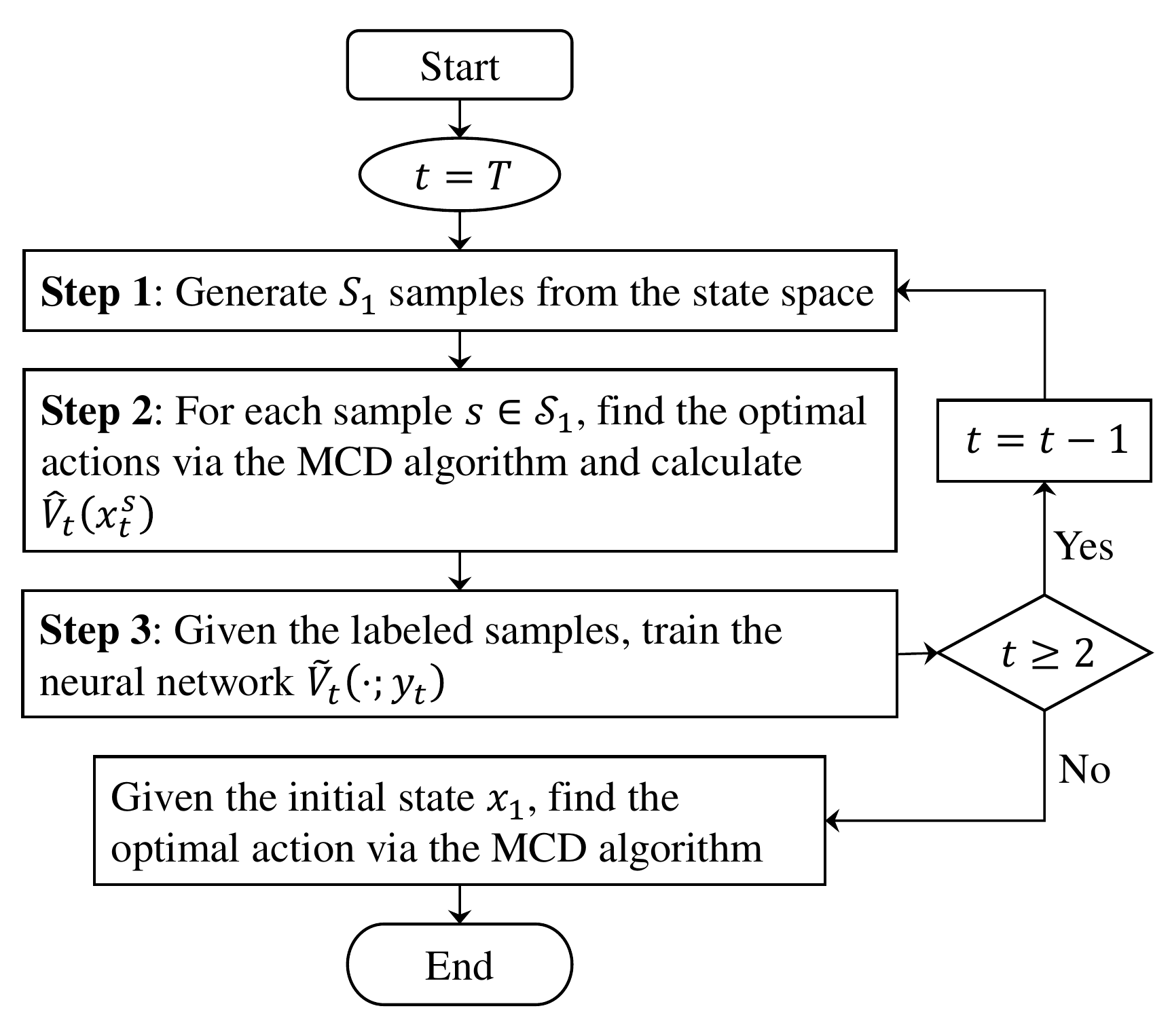}
			\par\end{centering}
		\caption{\label{fig:Alg_Flowchart}Flowchart of the NN-FVI algorithm combined with MCD}
	\end{figure}

	Cuts \eqref{eq:LLCuts_MP} are the mechanism which ensures convergence to a global optimum. However, they can theoretically cause slow convergence when the problem size is large. Yet, as we will see	in our numerical study, the MCD algorithm usually finds a high-performance solution (or even the optimal solution) in a relatively small number of iterations well before the gap between the lower/upper bound reaches the preset tolerance. Therefore, one can choose a suitable maximum number of iterations $\bar{m}$ to trade off between performance and computational budget. The flow of NN-FVI combined with MCD is summarized in Figure \ref{fig:Alg_Flowchart}, and the detailed procedure of MCD is presented in Algorithm \ref{alg:MCD_procedure}.
	
	{\small{}{}} 
	\begin{algorithm}[!tph]
		\caption{\label{alg:MCD_procedure}The MCD algorithm}
		
		\begin{algorithmic}[1] \Require{$\bar{m}, \epsilon$,	state $x_{t}$, initial solution $a_{t}^{1}$} \Ensure{$a_{t}^{*}$}
			\State{Initialize $m=1,V_{\text{lb}}=-\infty,V_{\text{ub}}=+\infty$}
			\While{$m<\bar{m}$ $\text{{\bf {and}}}$ $V_{\text{ub}}-V_{\text{lb}}<\epsilon$
			}
			\State{Evaluate $V^{m}=r_{t}\left(x_{t},a_{t}^{m}\right)+\gamma v_{t+1}(x_{t},a_{t}^{m};y_{t+1})+\gamma w_{0(t+1)}$ and set $V_{\text{lb}} \leftarrow \max\{V^{m},V_{\text{lb}}\}$ } 
			\State{Given
				$a_{t}^{m}$, construct Cuts \eqref{eq:LLCuts_MP} and \eqref{eq:validcut_all2} }
			\State{Solve
				Problem \eqref{eq:FP} and derive $a_{t}^{m+1}$ and $V_{\text{ub}}$}
			\State{$m\leftarrow m+1$} \EndWhile 
			\State{$a_{t}^{*}=\arg\max_{a\in\{a_{t}^{1},\ldots,a_{t}^{m}\}}\left\lbrace r_{t}(x_{t},a)+\gamma v_{t+1}(x_{t},a_{t}^{m};y_{t+1}) +\gamma w_{0(t+1)} \right\rbrace$}
		\end{algorithmic}
	\end{algorithm}
	
	\section{Application: Multi-facility Capacity Investment Problem}
	
	\label{sec:MCIP}
	
	We apply NN-FVI with MCD to solve a generic multi-facility
	capacity investment problem (MCIP) with discrete capacity decisions. Strategic
	capacity decisions are important to production firms because of the high
	expenditures entailed and the uncertainty about the future business
	environment (e.g. customer demand). To deal with uncertainty, it is wiser to
	adjust the capacity periodically as new information about the uncertainty is revealed, instead of establishing facilities
	with large capacities at the very beginning of the planning horizon. Real options analysis
	gives a framework to evaluate the value of a system
	with dynamically adjusted capacity. Namely, in each period,
	decision makers have the right (but not the obligation) to invest
	or salvage the capacity once new demand information is revealed \citep{dixit_investment_1994, huang_value_2009, cardin_strategic_2017, taghavi_lagrangian_2018}.
	
	To solve this problem, \cite{cardin_approach_2017} and \cite{zhao_decision_2018} have earlier proposed a decision rule-based multi-stage stochastic programming method, where the policy space of the problem is approximated by a family of if-then functions. The performance of the decision rule-based method is promising, but it only allows for capacity expansion decisions. We go in a different direction here and consider an MCIP where the capacity can be both expanded and contracted. We formulate the MCIP as a finite horizon MDP, and analyze the properties of its value functions. Then, we solve the problem via NN-FVI with MCD. Our comprehensive numerical experiments help to verify the performance of our algorithm.
	
	\subsection{Problem Description}
	
	We consider a multi-stage capacity investment problem with
	customers $\mathcal{I}\triangleq\left\{ 1,\ldots,I\right\} $ and
	facilities $\mathcal{\mathcal{N}}\triangleq\left\{ 1,\ldots,N\right\} $.
	In each period, customer demand is observed and then assigned to the various facilities to be satisfied. The objective is to maximize the net present value
	(NPV) by determining the optimal capacity investment plan and demand
	allocation plan over the planning horizon $\mathbb{T}$.
	Let $D_{it}$ denote the random demand from customer $i\in\mathcal{I}$
	in period $t\in\mathbb{T}$, and let $D_{t}=\left(D_{1t},\ldots,D_{It}\right)$ denote the vector of customer demands in period $t\in\mathbb{T}$.
	Let $d_{it}$ denote a realization of $D_{it}$, and the corresponding vector of realizations is $d_{t}=\left(d_{1t},\ldots,d_{It}\right)$. Define $\mathbb{D}\triangleq\bigcup_{t\in\mathbb{T}}\text{supp}\left(D_{t}\right)$ to be the support of the customer demand.
	
	Let $K_{0}=\left(K_{10},\ldots,K_{N0}\right)$
	be the initial capacity vector, and let $K_{t}=\left(K_{1t},\ldots,K_{Nt}\right)$ be the installed capacity
	vector at the end of period $t$. 
	Define $K^{\max}=\left(K_{1}^{\max},\ldots,K_{N}^{\max}\right)$
	to be a vector of capacity limits,
	so that the domain of the capacity is
	\[
	\mathbb{K}\triangleq\left\{ K\in\mathbb{Z}_{+}^{N}\left|0\le K\le K^{\text{max}}\right.\right\} .
	\]
	We summarize the notation for our model in Table \ref{tab:Notations_MCIP}.
	
	Our main assumptions on the MCIP are listed below. 
	\begin{assumption}
		\label{assu:Assumption1_Markov}The demand process $\left\{ D_{t}\right\}_{t\in\mathbb{T}}$
		is Markov. Without loss of generality,
		$D_{t}$ is independent of the capacity decisions $K_{t}$
		for all $t\in\mathbb{T}$, and $d_1$ is known.
	\end{assumption}
	
	Under Assumption \ref{assu:Assumption1_Markov}, we may let $P\left(\cdot\left|d_{t-1}\right.\right)$ denote the conditional probability distribution for demand so that $D_{t}\sim P\left(\cdot\left|d_{t-1}\right.\right)$
	for all $d_{t}\in\mathbb{D}$ and $t\in\mathbb{T}\backslash\{1\}$. If demand is continuous,
	we assume that $P\left(\cdot\left|d_{t-1}\right.\right)$ has a conditional probability density.
	
	\begin{assumption}
		\label{assu:Assumption2_Lip}The demand is non-negative and bounded
		(i.e., there exists $D^{\max}<\infty$ such that $0\le D_{t}\le D^{\max}$
		for all $t\in\mathbb{T}$). If the demand is continuous, the conditional probability densities $P\left(\cdot\left|d_{t-1}\right.\right)$ exist and
		are Lipschitz continuous for all $t\in\mathbb{T}\backslash\{1\}$. Specifically,
		there exists $L_{d}<\infty$ such that 
		\[
		\left|P\left(d_{t}\left|d_{t-1}\right.\right)-P\left(d_{t}\left|d_{t-1}^{\prime}\right.\right)\right|\le L_{d}\left\Vert d_{t-1}-d_{t-1}^{\prime}\right\Vert ,\ \ \forall d_{t-1},d_{t-1}^{\prime},d_{t}\in\mathbb{D},t\in\mathbb{T}\backslash\{1\}.
		\]
	\end{assumption}
	\noindent Assumption \ref{assu:Assumption2_Lip} is standard since real-world
	demand is always finite and their variation from one period to the
	next does not surge to infinity.
	
	\begin{assumption}
		The decision maker has the option to expand/contract the installed capacity
		from $K_{t-1}$ to $K_{t}$ at the end of each period $t\in\mathbb{T}$,
		and capacity adjustments are instantaneous.
	\end{assumption}
	
	We assume that the lead
	time for capacity adjustment (compared to the length of each period)
	is negligible. 
	
	\begin{assumption}
		\label{assu:adjustment_costs}The expansion cost and salvage value
		are linear with respect to the capacity, and the per unit expansion
		cost is not smaller than the per unit salvage value.
	\end{assumption}
	
	Though we consider an MCIP with linear expansion cost and
	salvage value, our proposed method can be extended to more general
	cases whose expansion cost/salvage value are nonlinear. In addition,
	the resale value of an asset is usually smaller than its purchase
	price, so the second statement of Assumption \ref{assu:adjustment_costs}
	is realistic. 
	
	\begin{table}
		\caption{\label{tab:Notations_MCIP}Notations for the multi-stage MCIP}
		\small
		\centering{}%
		\begin{tabular}{ll}
			\hline 
			\multicolumn{2}{c}{\emph{Indices and sets}}\tabularnewline
			{$i$} & {Index for customers}\tabularnewline
			{$n$} & {Index for facility}\tabularnewline
			{$t$} & {Index for period}\tabularnewline
			{$\mathcal{I}$} & {Set of customers, $i\in\mathcal{I}$, and $\left|\mathcal{I}\right|=I$}\tabularnewline
			{$\mathcal{N}$} & {Set of facilities, $n\in\mathcal{N}$, and $\left|\mathcal{N}\right|=N$}\tabularnewline
			{$\mathbb{T}$} & {Set of periods, $t\in\mathbb{T}$, and $\left|\mathbb{T}\right|=T$}\tabularnewline
			& \tabularnewline
			\multicolumn{2}{c}{\emph{Parameters}}\tabularnewline
			{$d_{it}$} & {Demand generated from customer $i$ in time $t$; the vector form is $d_{t}=\left(d_{1t},\ldots,d_{It}\right)$ }\tabularnewline
			{$K_{n0}$} & {Initial capacity of facility $n$; the vector form is $K_{0}=\left(K_{10},\ldots,K_{N0}\right)$}\tabularnewline
			{$\gamma$} & {Discount factor of time value of money, $0<\gamma<1$}\tabularnewline
			{$\hat{r}_{int}$} & {Unit revenue from satisfying customer $i$ with facility
				$n$ in time $t$}\tabularnewline
			{$b_{it}$} & {Unit penalty cost for unsatisfied customer $i$ in time $t$}\tabularnewline
			{$q_{nt}^{+}$} & {Coefficient parameters of per unit expansion cost of facility
				$n$ in time $t$}\tabularnewline
			{$q_{nt}^{-}$} & {Coefficient parameters of per unit salvage value of facility
				$n$ in time $t$}\tabularnewline
			& \tabularnewline
			\multicolumn{2}{c}{\emph{Variables}}\tabularnewline
			{$K_{nt}$} & {Capacity of facility $n$ in time $t\in\mathbb{T}\backslash\{1\}$; the vector form is
			$K_{t}=\left(K_{1t},\ldots,K_{Nt}\right)$}\tabularnewline
			{$z_{int}$} & {Amount of demand allocated from customer $i$ to facility $n$ in time $t$}\tabularnewline
			\hline 
		\end{tabular}
	\end{table}
	
	In each period, we allocate realized demand $d_{t}$ to the
	facilities, given the limits of the currently installed capacity
	$K_{t-1}$. A unit penalty cost $b_{it}$ is incurred for all unsatisfied demand
	from customer $i$. Let $z_{int}$ denote the amount
	of demand from customer $i$ allocated to facility $n$ in period
	$t$, and let $\hat{r}_{int}$ be the corresponding revenue. The operating
	profit $Q_{t}\left(K_{t-1},d_{t}\right)$ in period $t\in\mathbb{T}$
	is given by the following linear program (LP):
	\begin{subequations}
		\begin{alignat}{2}
		Q_{t}\left(K_{t-1},d_{t}\right) \triangleq\max_z & \sum_{i\in\mathcal{I}}\sum_{n\in\mathcal{N}}\hat{r}_{int}z_{int}-\sum_{i\in\mathcal{I}}b_{it}\left(d_{it}-\sum_{n\in\mathcal{\mathcal{N}}}z_{int}\right)\label{eq:Lp_obj} & \\
		\text{s.t.}\ \  & \sum_{i\in\mathcal{I}}z_{int}\le K_{n\left(t-1\right)}, & \forall n\in\mathcal{N},\label{eq:Cap_con}\\
		& \sum_{n\in\mathcal{N}}z_{int}\le d_{it}, & \forall i\in\mathcal{I},\label{eq:D_con}\\
		& z_{int}\ge0, & \forall i\in\mathcal{I},n\in\mathcal{N}.
		\end{alignat}
	\end{subequations}
	The	objective in \eqref{eq:Lp_obj} is to maximize the current rewards, which consist of the revenue minus the penalty for unsatisfied demand. Constraints \eqref{eq:Cap_con} and \eqref{eq:D_con} are the capacity and demand constraints, respectively. Note that the allocation decisions $z_{int}$ depend on the current state $\left(K_{t-1},d_{t}\right)$ only. 
	
	Let $q_{nt}^{+}$ and $q_{nt}^{-}$ denote the unit expansion cost
	and unit salvage value for facility $n\in\mathcal{N}$ respectively.
	The reward is then
	\begin{equation}
	r_{t}\left(\left(K_{t-1},d_{t}\right),K_{t}\right)\triangleq
	Q_{t}\left(\left(K_{t-1},d_{t}\right)\right) - \sum_{n\in\mathcal{N}}\max\left\{ q_{nt}^{-}\left(K_{nt}-K_{n\left(t-1\right)}\right),q_{nt}^{+}\left(K_{nt}-K_{n\left(t-1\right)}\right)\right\}.\label{eq:cost_function}
	\end{equation}
	According to Assumption \ref{assu:adjustment_costs}, we have $q_{nt}^{+}\ge q_{nt}^{-}$
	for all $n\in\mathcal{N}$ and $t\in\mathbb{T}$.
	
	This MCIP can be modeled as
	an MDP where the state in period $t\in\mathbb{T}$ is $x_{t}\triangleq\left(K_{t-1},d_{t}\right)$ (i.e., the installed
	capacity in time $t-1$ and the realized demands) and the action
	is the adjusted capacity $K_{t}$. The state space of our problem is therefore $\mathbb{X}\triangleq\left\{ \left(K,d\right)\in\mathbb{K}\times\mathbb{D}\right\} ,$
	and the action space is $\mathbb{A} \triangleq \mathbb{K}$ for all $t\in\mathbb{T}$.
	Without loss of generality, we assume the initial state for the MCIP
	$x_1 = \left(K_{0},d_{1}\right)$ is known, and the system salvages
	all installed capacity at the end of period $T$ so $K_{T}\equiv\boldsymbol{0}_{N}$.
	The DP equations for MCIP are then:
	\begin{subequations}\label{eq:MCIP_value_function_1}
		\begin{align}
		V_{T}\left(K_{T-1},d_{T}\right)= & \max\limits _{K_{t}\in\mathbb{K}}\left\{ r_{T}\left(\left(K_{T-1},d_{T}\right),K_{T}\right)\right\},\\
		V_{t}\left(K_{t-1},d_{t}\right)= & \max\limits _{K_{t}\in\mathbb{K}}\left\{ r_{t}\left(\left(K_{t-1},d_{t}\right),K_{t}\right)+\gamma\mathbb{E}\left[V_{t+1}\left(K_{t},D_{t+1}\right)\left|d_{t}\right.\right]\right\},\,\forall t\in\mathbb{T}\backslash\left\{ T\right\}.
		\end{align}
	\end{subequations}
	If the demand has finite support, then complexity
	of VI for Eq. \eqref{eq:MCIP_value_function_1} is $O\left(\left|\mathbb{X}\right|^{2}\times\left|\mathbb{K}\right|\times T\right)$,
	where $\left|\mathbb{X}\right|$ has dimension $I+N$ and $\mathbb{K}$
	has dimension $N$. For example, consider a system with four customers, three facilities, and planning horizon $T=10$. If $K^{\max}=9$, and customer demand is integer-valued ranging from $1$ to $10$, then
	the complexity of VI is $\left(10^{3}\times10^{4}\right)^{2}\times10^{3}\times10=10^{18}$.
	This observation suggests that VI is intractable even for a medium size MCIP problem. If demand is continuous, then the state space is infinite and VI cannot even
	be done exactly. 
	
	\subsection{Problem Analysis}
	
	In this subsection, we analyze the value functions of the MCIP corresponding to the DP Eq. \eqref{eq:MCIP_value_function_1}. Then, we check if NN-FVI can solve the problem properly. We first show that the assumptions in Section \ref{sec:Preliminaries} are satisfied.
	\begin{proposition}
		\label{prop:check_assumptions}Problem \eqref{eq:MCIP_value_function_1}
		satisfies Assumptions \ref{assu:MDP_finite_and_discounted}--\ref{assu:MDP_lip_transition}. 
	\end{proposition}
	
	We observe that the value functions $V_{t}\left(\cdot\right)$ for all $t\in\mathbb{T}$ provided by Eq. \eqref{eq:MCIP_value_function_1} are defined over $\mathbb{X}$. However, $\mathbb{X}$ is not a connected set since $\mathbb{K}$ is finite. In addition, if the demand is discrete, then this aspect of the state space is finite and discrete as well. We wish to extend $\mathbb{K}$ to its smallest connected superset on which we will construct a set of extended value functions which are easier to analyze.
	
	We first define
	\[
	\bar{\mathbb{X}}\triangleq \bar{\mathbb{K}}\times\mathbb{D} \ \ \text{where}\ \ \bar{\mathbb{K}}\triangleq\left\{ K\in\mathbb{R}^{N}\left|0\le K\le K^{\text{max}}\right.\right\}.
	\]
	It is automatic that $\mathbb{X}\subset\bar{\mathbb{X}}$. Note that $\mathbb{D}$
	can be either continuous or discrete in the above display. Now,
	we define extended value functions $\bar{V}_{t}:\bar{\mathbb{X}}\mapsto\mathbb{R}$
	for all $t\in\mathbb{T}$. If $\bar{V}_{t}\left(\cdot\right)$ can
	recover the exact values of $V_{t}\left(\cdot\right)$ on $\mathbb{X}$,
	then we can analyze $\bar{V}_{t}\left(\cdot\right)$
	to learn about $V_{t}\left(\cdot\right)$. The analysis of $\bar{V}_{t}\left(\cdot\right)$
	is easier because $\bar{\mathbb{K}}$ is a connected set.
	
	The DP equations for the extended value functions
	at $\left(K_{t-1},d_{t}\right)\in\bar{\mathbb{X}}$ are:
	\begin{subequations}\label{eq:ValueFunction-Appro}
		\begin{align}
		\bar{V}_{T}\left(K_{T-1},d_{T}\right)\triangleq & \max\limits _{K_{T}\in\mathbb{K}}\left\{ r_{T}\left(\left(K_{T-1},d_{T}\right),K_{T}\right)\right\},\\
		\bar{V}_{t}\left(K_{t-1},d_{t}\right)\triangleq & \max\limits _{K_{t}\in\mathbb{K}}\left\{ r_{t}\left(\left(K_{t-1},d_{t}\right),K_{t}\right)+\gamma\mathbb{E}\left[\left.\bar{V}_{t+1}\left(K_{t},D_{t+1}\right)\right|d_{t}\right]\right\} ,\, \forall t\in\mathbb{T} \backslash \left\{ T \right\}.
		\end{align}
	\end{subequations}
	Note that the values of $r_{t}\left(\cdot\right)$ are well-defined for all $\left(K_{t-1},d_{t}\right)\in\bar{\mathbb{X}}$.
	
	Next we show that the value functions $V_{t}\left(\cdot\right)$ can be recovered from the extended functions $\bar{V}_{t}\left(\cdot\right)$ on $\mathbb{X}$.
	\begin{proposition}
		\label{prop:equivalence_value_functions}$\bar{V}_{t}\left(K_{t-1},d_{t}\right)=V_{t}\left(K_{t-1},d_{t}\right)$
		for all $\left(K_{t-1},d_{t}\right)\in\mathbb{X}$ and $t\in\mathbb{T}$.
	\end{proposition}
	
	We analyze the structure of the extended value functions $\left\{ \bar{V}_{t}\left(\cdot\right)\right\}_{t\in\mathbb{T}}$ by backward induction.
	
	\renewcommand{\labelenumi}{(\roman{enumi})}
	\begin{proposition}
		\label{prop:Piecewise Value Function}Let $\left\{ \bar{V}_{t}\left(\cdot\right)\right\}_{t\in\mathbb{T}}$ be the extended value functions defined in Eq. \eqref{eq:ValueFunction-Appro}.
		\begin{enumerate}
			\item $\bar{V}_{T}\left(\cdot\right)$ is piecewise linear.
			\item Suppose $\mathbb{D}$ is finite, then $\bar{V}_{t}\left(\cdot\right)$ is piecewise linear for all $t\in\mathbb{T}\backslash\left\{ T\right\}$.
		\end{enumerate}
	\end{proposition}
	
	\begin{remark}
		\label{rem:non_concave_V}Note that $\bar{V}_{t}\left(\cdot,d_{t}\right)$
		for all $t\in\mathbb{T}\backslash\left\{ T\right\} $ is not necessarily concave
		in $K_{t-1}\in\bar{\mathbb{K}}$ given $d_{t}\in\mathbb{D}$.
		As $\bar{V}_{T}\left(\cdot,d_{T}\right)$ is concave in $K_{T-1}\in\bar{\mathbb{K}}$
		given $d_{T}\in\mathbb{D}$, $\mathbb{E}\left[\bar{V}_{T}\left(\cdot,D_{T}\right)\left|d_{T-1}\right.\right]$
		is concave in $K_{T-1}\in\bar{\mathbb{K}}$ given $d_{T-1}$ (see
		the proof in Proposition \ref{prop:Piecewise Value Function}). Yet, $\bar{V}_{T-1}\left(\cdot,d_{T-1}\right)$
		may fail to be concave as it is a finite maximum of concave functions.
	\end{remark}
	
	According to Proposition \ref{prop:Piecewise Value Function}(i),
	the extended value function $\bar{V}_{T}\left(\cdot\right)$ is piecewise linear if the capacity is defined over a connected space $\bar{\mathbb{K}}$.
	However, for $t\in\mathbb{T}\backslash\left\{ T\right\} $, Proposition
	\ref{prop:Piecewise Value Function}(ii) the value
	functions $\bar{V}_{t}\left(\cdot\right)$ are only piecewise linear	when the demand is discrete. Fortunately, if the demand is continuous,
	we can use Monte-Carlo simulation to generate finitely many state transitions to approximate the expectation in Eq. \eqref{eq:MCIP_value_function_1}.
	In this setting, the approximate extended value functions will be piecewise
	linear. However, as observed in Remark \ref{rem:non_concave_V}, we may need
	to solve a non-convex optimization problem for each $t\in\mathbb{T}\backslash\left\{ T\right\} $.
	
	\subsection{\label{sec:Numerical-Studies}Numerical Studies}
	
	We now test the numerical performance of NN-FVI in solving MCIP in three parts. First, we compare the performance of the ReLU-based NN-FVI to NNs with other activation functions. This part is intended to verify that the value functions of MCIP have some piecewise linear structure, so that ReLU outperforms other activation functions. Second, we compare the performance of the proposed MCD algorithm with exhaustive enumeration (i.e., the brute-force method) and the integer L-shaped algorithm in the action selection problem in NN-FVI. Third, we combine NN-FVI with MCD and test its performance in a case study where we analyze its economic performance over an inflexible counterpart. The inflexible counterpart does not have the option to change the capacity dynamically. Instead, it is modeled as a two-stage capacity investment problem and it is solved with Benders decomposition (see \cite{zhao_decision_2018} for implementation details). These numerical studies are performed on a workstation with an Intel Xeon 5218 processor and 32 GB RAM in the Matlab R2018a environment. The NNs are trained by the Levenberg--Marquardt algorithm via the NN toolbox of Matlab \citep{beale_neural_2018}. 
	
	\subsubsection{ReLU Outperforms Other Activation Functions in Solving MCIP}
	
	In this subsection, we test the performance of NN-FVI with different
	types of activation functions, including ReLU, tanH (hyperbolic function),
	and sigmoid. Here, the action selection problem is solved by the brute-force
	method as MCD is not applicable for NNs using tanH and sigmoid
	functions. 
	
	In Table \ref{tab:ActFun}, we test a small-scale case (Case 1.1)
	with discrete demands that it is solvable by DP, and compare the approximate
	objective values obtained from NN-FVI to the exact objective values
	obtained by DP. As can be seen, the exact ENPV computed by DP is 357.9, and the approximate objective computed by NN-FVI
	with ReLU is 357.7. However, the approximate ENPVs computed by NN-FVI
	with tanH and sigmoid activation functions are 441.3 and 628.9 respectively,
	both of which are far from the exact value. \vspace{-0.05cm}
	
	\begin{table}[H]
		\caption{\label{tab:ActFun}Comparisons of NN-FVI with different activation functions for small-scale cases}
		
		\begin{centering}
			\begin{tabular}{cccccc}
				\toprule 
				& Algorithms & Activation fun. & CPU time & Obj. values{*} & Relative gaps\tabularnewline
				\midrule 
				\multirow{4}{*}{$\begin{array}{c}
					\text{Case 1.1}\\
					\left(I,N,T\right)=\left(2,2,2\right)
					\end{array}$} & DP & - & 160 s & 357.9 & -\tabularnewline
				& NN-FVI & ReLU & 6 s & 357.7 & \textless 0.1\%\tabularnewline
				& NN-FVI & tanH & 10 s & 441.3 & 23.3\%\tabularnewline
				& NN-FVI & Sigmoid & 10 s & 628.9 & 75.7\%\tabularnewline
				\bottomrule
			\end{tabular}
		\end{centering}
		\raggedright{}{\small{}{*}The ENPVs are computed from the approximated objective values directly.}
	\end{table}	
	For a larger scale case, DP is not applicable due to the curse
	of dimensionality. Instead, we benchmark using an inflexible two-stage stochastic capacity
	investment model. We perform out-of-sample
	tests on the optimal policies computed by both the inflexible method
	and NN-FVI, on an identical sample set with 10,000 sample paths; in
	this case, a better policy should lead to a higher ENPV in the out-of-sample
	tests. In Table \ref{tab:ActFun-1}, simulation results for a medium-scale
	case with $\left(I,N,T\right)=\left(4,3,11\right)$ indicate that
	the optimal policy computed by NN-FVI with ReLU outperforms the policies
	computed by NNs with tanH and sigmoid activation functions. Also, the CPU
	time for NN-FVI with ReLU is much less than NN-FVI with
	tanH and sigmoid. We suspect this is because the piecewise linear activation functions are
	easier to use in computation compared to the nonlinear ones. We can
	conclude that ReLU leads to more accurate NNs compared to tanH and sigmoid.
	
	\begin{table}[H]
		\caption{\label{tab:ActFun-1}Comparisons of NN-FVI with different activation functions for medium-scale cases}		
		\centering{}%
		\begin{tabular}{cccccc}
			\toprule 
			& Algorithms & Activation fun. & CPU time & ENPV & VoF\tabularnewline
			\midrule
			\multirow{4}{*}{$\begin{array}{c}
				\text{Case 1.2}\\
				\left(I,N,T\right)=\left(4,3,11\right)
				\end{array}$} & Inflexible design & - & 306 s & 1530.9 & -\tabularnewline
			& NN-FVI & ReLU & 21 040 s & 1634.9 & 104.0\tabularnewline
			& NN-FVI & tanH & 25 624 s & 1625.4 & 94.5\tabularnewline
			& NN-FVI & Sigmoid & 27 452 s & 1518.2 & $-12.7$\tabularnewline
			\bottomrule
		\end{tabular}
	\end{table}

	\subsubsection{The Action Selection Procedure Can be Solved in Reasonable Time}
	
	In this subsection, we compare the performance of our MCD algorithm
	with two alternatives---(1) the brute-force algorithm and (2)
	the integer L-shaped algorithm. In the integer L-shaped algorithm,
	only Cuts \eqref{eq:LLCuts_MP} (integer optimality cuts) are added
	in each iteration and our valid cuts are not included. 
	
	Seven NNs with different sizes are
	randomly generated to create Cases 2.1--2.7, each is formulated as Problem \eqref{eq:OAS_obj},
	and then solved by the aforementioned methods. The performance of the
	algorithms is measured in terms of the CPU time and their best-found
	objective value achieved before the algorithm is terminated. As can
	be seen, in Cases 2.1 and 2.2, the CPU time of the brute-force method
	is around 11 seconds. Both the integer L-shaped algorithm and MCD converge to the optimal solution in around 6-7 seconds, which is slightly faster than the brute-force method. This is not surprising since the brute-force method is plausible when the problem size is small. By increasing
	the number of facilities from three to five, we can see from Table
	\ref{tab:Sim_result-MCD} that the CPU time of the brute-force method
	increases significantly. On the contrary, MCD achieves the global
	optimal/near-optimal solutions within 10 seconds. Based on this evidence, we speculate that	the CPU time of the MCD algorithm might be much less than the brute-force method in even larger scale cases. 
	
	\begin{table}
		\caption{\label{tab:Sim_result-MCD}Simulation results of the MCD algorithm}
		
		\begin{centering}
			\small
			\begin{tabular}{ccccccc}
				\toprule 
				& {Algorithms} & {Stop criterion{*}} & {Iterations} & {CPU time} & {Exp. profits} & {Relative gaps}\tabularnewline
				
				\midrule 
				{Case 2.1} & {Brute-force} & {-} & {-} & {11.3 s} & {2506.0} & {-}\tabularnewline
				{($N=3$)} & {Integer L-shaped} & {0.35\%/100 steps} & {100} & {6.7 s} & {2506.0} & {0\%}\tabularnewline
				& {MCD} & {0.35\%/100 steps} & {100} & {7.8 s} & {2506.0} & {0\%}\tabularnewline
				
				\midrule
				{Case 2.2} & {Brute-force} & {-} & {-} & {11.2 s} & {6731.9} & {-}\tabularnewline
				{($N=3$)} & {Integer L-shaped} & {0.35\%/100 steps} & {100} & {6.5 s} & {6721.6} & {0.15\%}\tabularnewline
				& {MCD} & {0.35\%/100 steps} & {100} & {7.6 s} & {6731.9} & {0\%}\tabularnewline
				
				\midrule
				{Case 2.3} & {Brute-force} & {-} & {-} & {47.9 s} & {$119.3\times10^{6}$} & {-}\tabularnewline
				{($N=4$)} & {Integer L-shaped} & {0.35\%/200 steps} & {200} & {18.7 s} & {$79.8\times10^{6}$} & {33.11\%}\tabularnewline
				& {MCD} & {0.35\%/100 steps} & {100} & {8.5 s} & {$119.2\times10^{6}$} & {\textless 0.1\%}\tabularnewline
				
				\midrule
				{Case 2.4} & {Brute-force} & {-} & {-} & {48.6 s} & {$466.2\times10^{6}$} & {-}\tabularnewline
				{($N=4$)} & {Integer L-shaped} & {0.35\%/100 steps} & {100} & {14.7 s} & {$465.0\times10^{6}$} & {0.26\%}\tabularnewline
				& {MCD} & {0.35\%/100 steps} & {100} & {8.3 s} & {$465.6\times10^{6}$} & {0.13\%}\tabularnewline
				
				\midrule
				{Case 2.5} & {Brute-force} & {-} & {-} & {643.4 s} & {$256.8\times10^{6}$} & {-}\tabularnewline
				{($N=5$)} & {Integer L-shaped} & {0.35\%/100 steps} & {100} & {14.3 s} & {$233.1\times10^{6}$} & {7.30\%}\tabularnewline
				& {MCD} & {0.35\%/100 steps} & {100} & {11.0 s} & {$256.8\times10^{6}$} & {0\%}\tabularnewline
				& {MCD} & {0.35\%/50 steps} & {50} & {4.3 s} & {$256.8\times10^{6}$} & {0\%}\tabularnewline
				
				\midrule  
				{Case 2.6} & {Brute-force} & {-} & {-} & {764.7 s} & {$113.8\times10^{6}$} & {-}\tabularnewline
				{($N=5$)} & {Integer L-shaped} & {0.35\%/100 steps} & {100} & {55.2 s} & {$112.0\times10^{6}$} & {1.54\%}\tabularnewline
				& {MCD} & {0.35\%/100 steps} & {40} & {4.4 s} & {$113.8\times10^{6}$} & {0\%}\tabularnewline
				
				\midrule
				{Case 2.7} 	& {Brute-force} & {-} & {-} & {815.5 s} & {$887.6 \times10^{6}$} & {-}\tabularnewline
				{($N=5$)} & {Integer L-shaped} & {0.35\%/100 steps} & {100} & {66.2 s} & {$886.7\times10^{6}$} & {\textless 0.1\%}\tabularnewline
				& {MCD} & {0.35\%/100 steps} & {14} & {1.6 s} & {$887.1\times10^{6}$} & {\textless 0.1\%}\tabularnewline
				\bottomrule
			\end{tabular}
			\par\end{centering}
		\raggedright{}{\footnotesize{}{*}The algorithm is stopped if the relative gap between the lower bound and the upper bound is smaller than 0.35\% }\emph{\footnotesize{}or}{\footnotesize{}	the total number of iteration reaches 100.}{\footnotesize\par}
	\end{table}
	
	We now compare the convergence of MCD and the integer L-shaped algorithm.
	As can be seen, MCD converges to the global optimum
	in fewer iterations than the integer L-shaped algorithm. In Case
	2.1, where the number of facilities is $N=3$, the integer
	L-shaped method can converge to the global optimum within 100
	iterations. However, when the number of facilities increases, the
	best-found solutions by the integer L-shaped method within 100 iterations
	are sub-optimal. In particular, the relative gaps of the best-found
	objective values and the global optima for Case 2.3 ($N=4$), Case 2.5 ($N=5$), and Case 2.6 ($N=5$)	are $33.11\%$, $7.30\%$, and $1.54\%$ respectively. In contrast, MCD can find
	the global optimum (or a near-optimal solution) within 100 steps.
	These numerical results verify that MCD ensures faster convergence by using the
	valid cuts with gradient information (i.e., Cuts \eqref{eq:validcut_all2}).
	
	\begin{figure}

		\begin{centering}
			\includegraphics[scale=0.50]{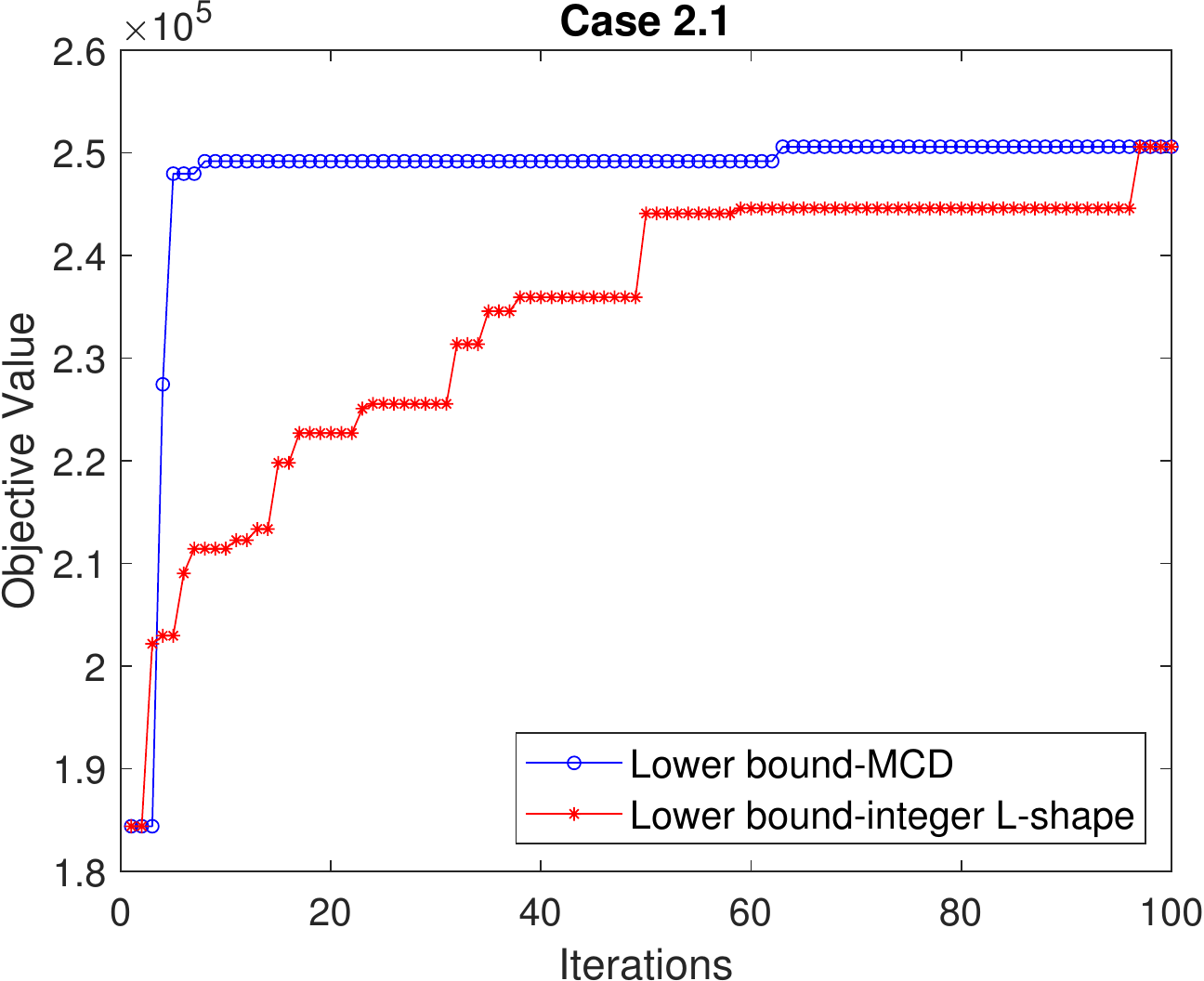}
			\includegraphics[scale=0.50]{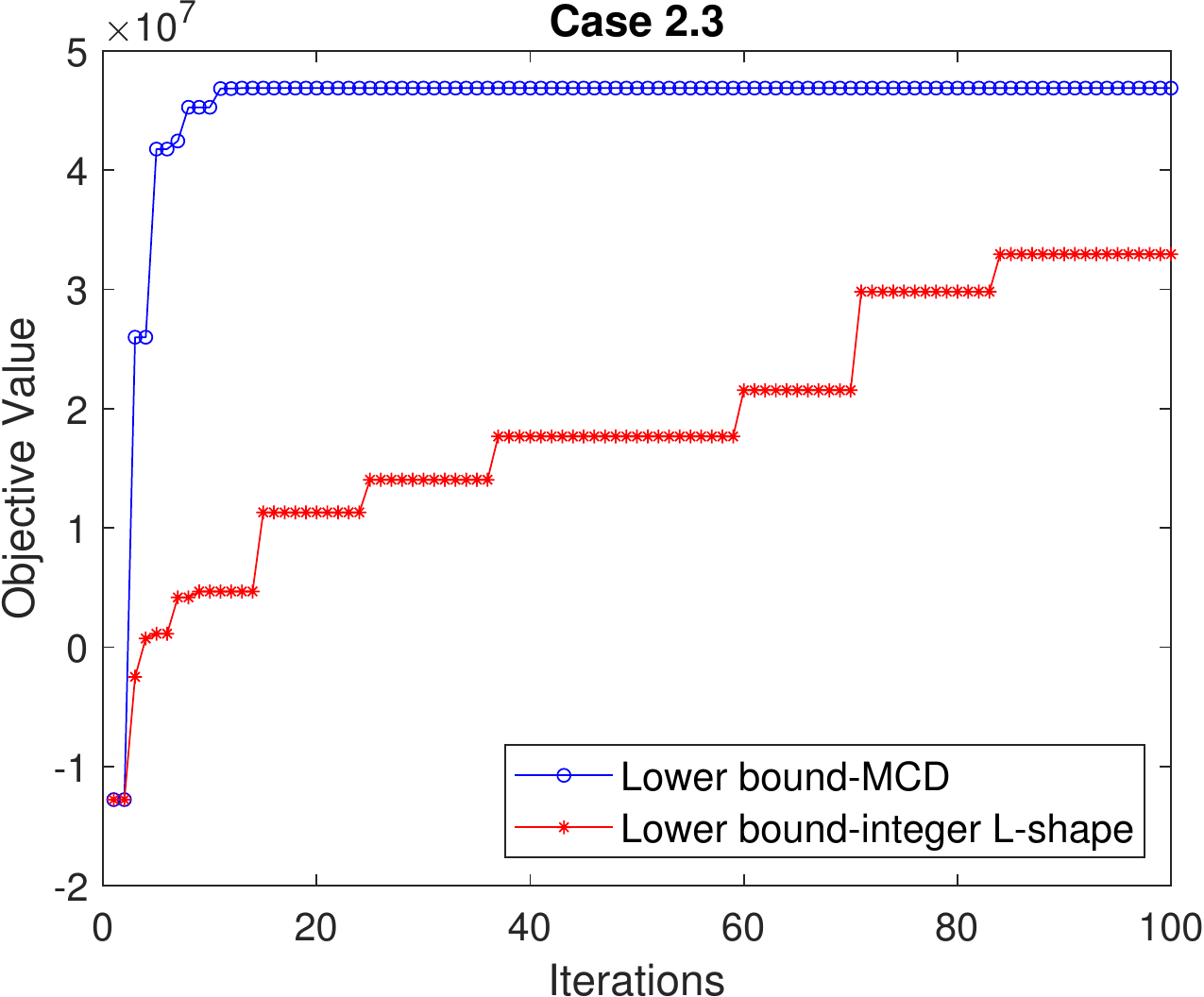}
			\par\end{centering}
		\vspace{0.5cm}
		\begin{centering}
			\includegraphics[scale=0.50]{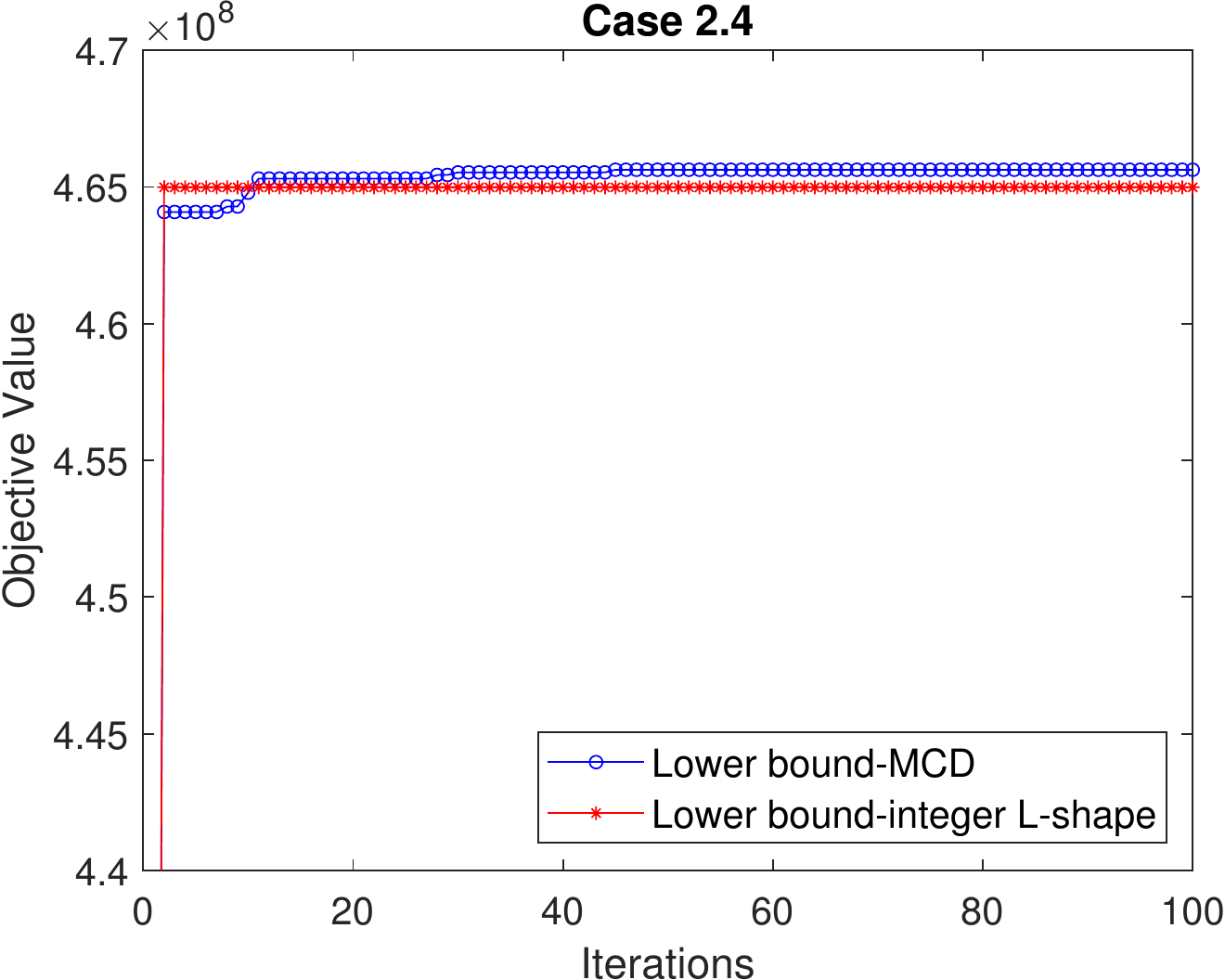}
			\includegraphics[scale=0.50]{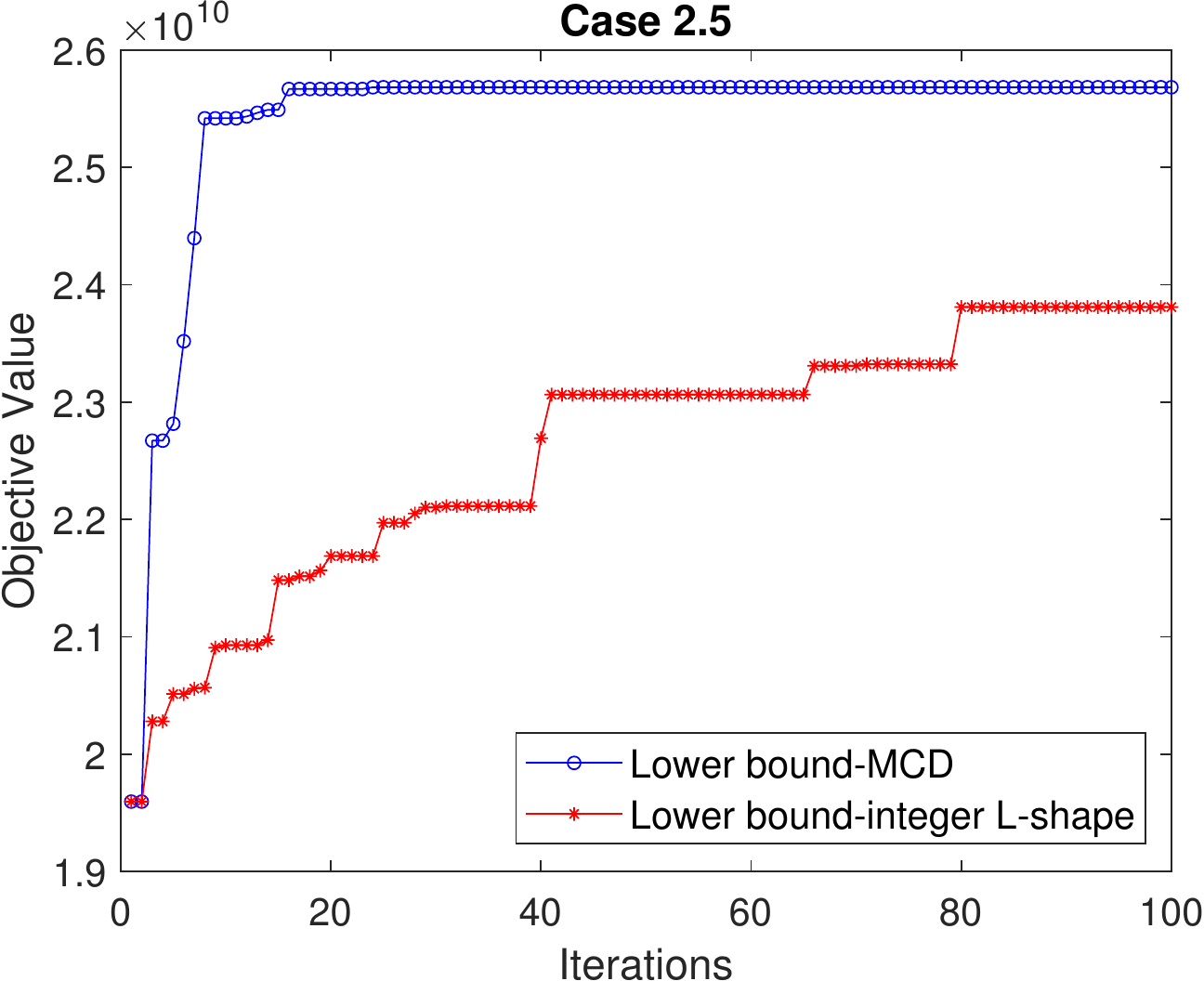}
			\par\end{centering}
		\vspace{0.5cm}
		\begin{centering}
			\includegraphics[scale=0.50]{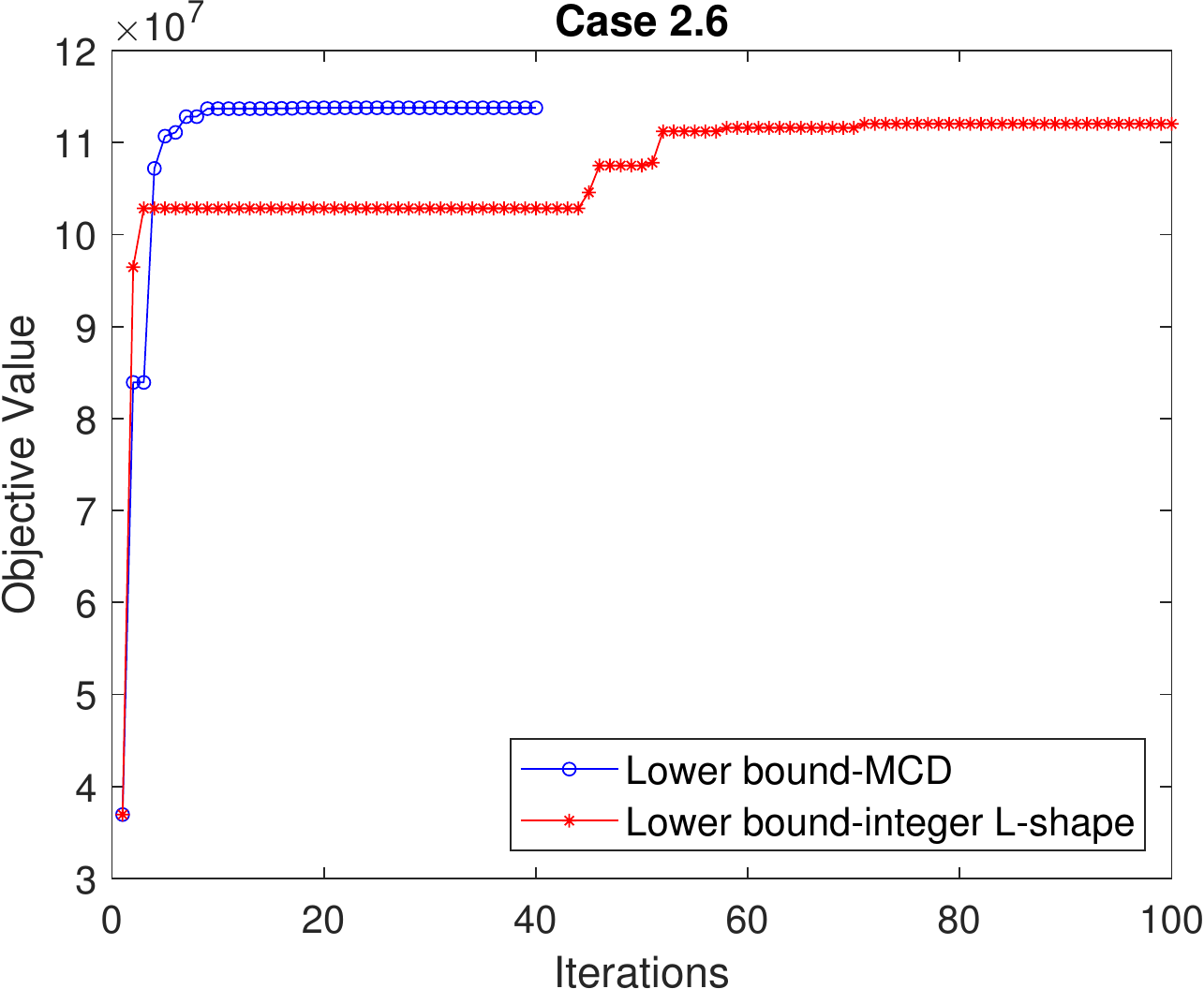}
			\includegraphics[scale=0.50]{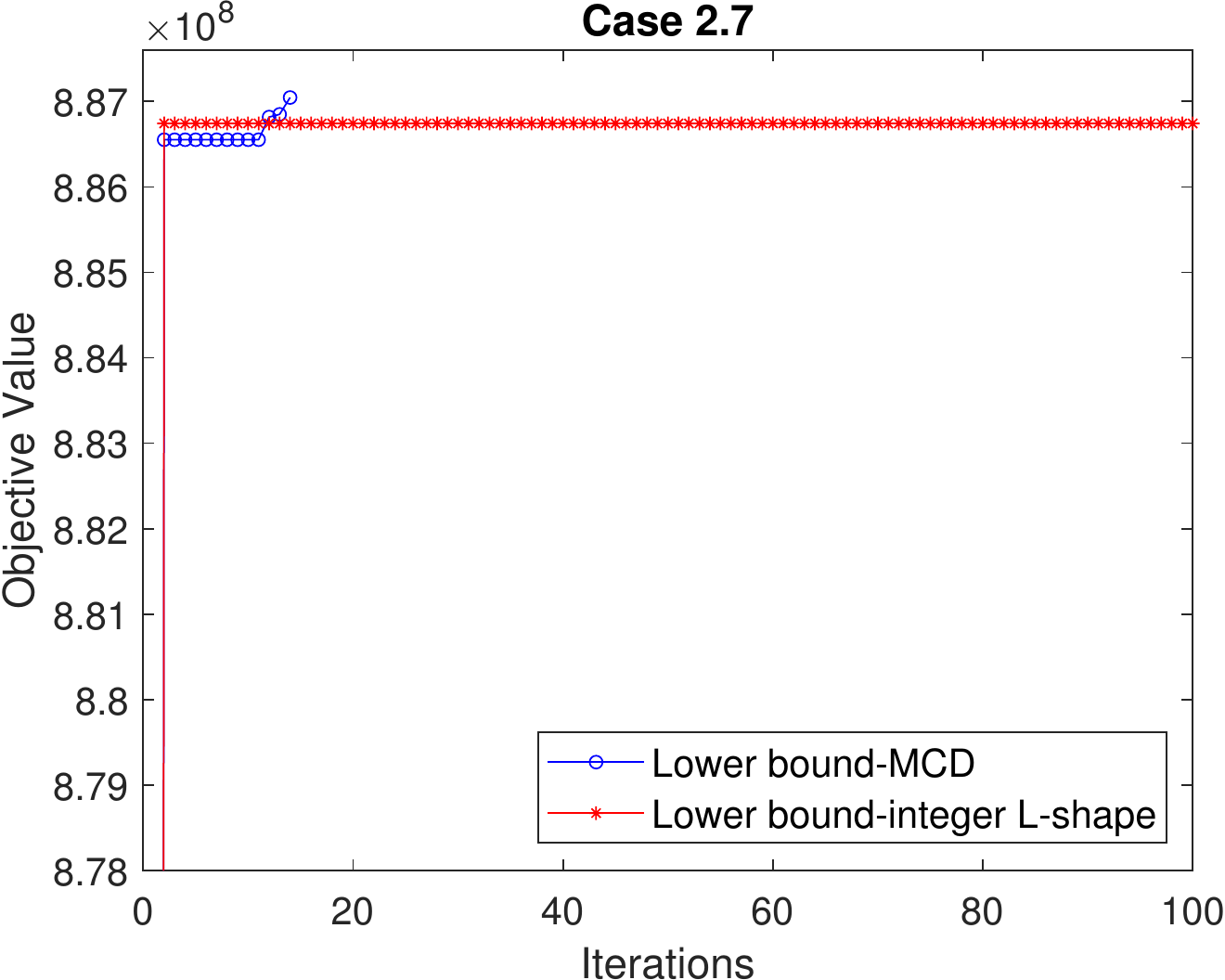}
			\par\end{centering}
		\caption{Convergence analyses of the MCD algorithm}
	\end{figure}

	\subsubsection{Verify the Performance of NN-FVI with MCD in a Real Case}
	
	This case study on a multi-facility waste-to-energy (WTE) system is
	adapted from \citep{zhao_decision_2018}. In \citep{zhao_decision_2018}, the
	system can only expand capacity, while in this
	paper we also allow capacity to contract. In addition, we adjust the
	original case into a smaller one where we can perform practical sensitivity
	analysis. The WTE system has four candidate sites located in different
	sectors. The facilities at each site are able to dispose of food waste
	collected from each sector by using an anaerobic digestion technique,
	which transforms the food waste into electricity. Undisposed waste
	will be subjected to further treatment via landfill, incurring greater
	disposal costs and penalties. The revenue comes from selling the
	electricity and the salvage value of the contracted capacity. The costs consist of the disposal, penalty, and transportation costs, as well as capacity expansion costs and land rental fees (incurred once the capacity is non-zero). The data and parameters of
	the case study are found in \citep{zhao_decision_2018}. 
	
	Sensitivity analysis is implemented on the ratio of the per unit salvage value with the per unit expansion cost (S/E ratio) $q_{nt}^{-}/q_{nt}^{+}$. If the S/E ratio is one, it means that the salvage value for per unit capacity is equivalent to the per unit expansion cost. As can be seen in Table \ref{tab:SensitivityAnalysis_S/E}, the flexible design outperforms the inflexible one under most parameter settings. When $\gamma = 0.862$ and the S/E ratio is 0, the expected net present value (ENPV) of the inflexible design is 23.9 million, while the ENPV of the flexible design is 35.5 million. In this case, flexibility improves the system performance by 48.5\%. However, the percentage improvement	of the system with a flexible design over an inflexible design decreases as the S/E ratio increases. The improvement decreases to 20.3\% when the S/E ratio increases to 0.99. On the other hand, we see that if the discount factor $\gamma$ and the S/E ratio are both close to one, then the improvement is only $1.0 \%$. This is reasonable since in this case the decision maker can establish facilities with large capacity in the beginning and then salvage all of them in the last period, without suffering a significant loss. In other words, the inflexible model may yield a high-performance solution for the problem when both $\gamma$ and the S/E ratio are close to one. 
	
	\begin{table}[H]
		\caption{\label{tab:SensitivityAnalysis_S/E}Sensitivity analysis given different
			S/E ratio and gamma}
		\small
		\centering{}%
		\begin{tabular}{cccccc}
			\toprule 
			\multirow{2}{*}{$\begin{array}{c}
				\text{Gamma}\\
				\left(\gamma\right)
				\end{array}$} & \multirow{2}{*}{$\begin{array}{c}
				\text{S/E ratio}\\
				\left(q_{nt}^{-}/q_{nt}^{+}\right)
				\end{array}$} & \multicolumn{2}{c}{ENPV $\left(\times10^{6}\text{\ S\$}\right)$} & \multirow{2}{*}{$\begin{array}{c}
				\text{VoF}\\
				\left(\times10^{6}\text{\ S\$}\right)
				\end{array}$} & \multirow{2}{*}{$\begin{array}{c}
				\text{Improvement}\\
				\left(\%\right)
				\end{array}$}\tabularnewline
			\cmidrule{3-4} \cmidrule{4-4} 
			&  & Inflexible design & NN-FVI &  & \tabularnewline
			\midrule 
			\multirow{5}{*}{0.862} 
			& 0 	& 23.9 & 35.5 & 11.6 & 48.5\tabularnewline
			& 0.25 	& 33.7 & 45.7 & 12.0 & 35.6\tabularnewline
			& 0.50 	& 43.8 & 58.0 & 14.2 & 32.4\tabularnewline
			& 0.75 	& 54.3 & 66.9 & 12.6 & 23.2\tabularnewline
			& 0.99 	& 64.9 & 78.1 & 13.2 & 20.3\tabularnewline
			\midrule
			\multirow{5}{*}{0.923} 
			& 0 	& 93.8  & 110.0 & 16.2 & 17.3\tabularnewline
			& 0.25	& 115.1 & 131.1 & 16.0 & 13.9\tabularnewline
			& 0.50	& 137.2 & 156.8 & 19.6 & 14.3\tabularnewline
			& 0.75	& 160.4 & 176.0 & 15.6 & 9.7\tabularnewline
			& 0.99	& 184.0 & 194.8 & 10.8 & 5.9\tabularnewline
			\midrule
			\multirow{5}{*}{0.99} 
			& 0		& 204.3 & 224.4 & 20.1 & 9.8\tabularnewline
			& 0.25 	& 248.9 & 268.0 & 19.1 & 7.7\tabularnewline
			& 0.50 	& 295.2 & 312.3 & 17.1 & 5.8\tabularnewline
			& 0.75 	& 345.9 & 361.1 & 15.2 & 4.4\tabularnewline
			& 0.99 	& 398.4 & 402.4 & 4.0  & 1.0\tabularnewline
			\bottomrule
		\end{tabular}
	\end{table}

	\subsection{\label{sec:Extensions}Extensions of NN-FVI with MCD}
	
	We focus on solving MCIP in this paper, but our method is applicable
	to many other problems where the action space is finite and high-dimensional. We mention some variations where our method can still be used.
	
	\renewcommand{\labelenumi}{\arabic{enumi}.}
	\begin{enumerate}
		\item \emph{Lead time}. In MCIPs where the capacity adjustment
		has nonzero lead time, we can express capacity under construction as
		part of the state variables. In this case, we may need additional inputs for our two-layer
		NN, but the action selection problem
		is unaffected.
		\item \emph{Uncertain rewards/costs}. The proposed method can also solve
		MCIPs with uncertain parameters. For example, if the rewards are uncertain, then we can model them as state variables.
		\item \emph{Different types of capacity adjustment costs}. If the capacity
		adjustment costs contain fixed costs or are non-convex, the proposed
		MCD algorithm can still solve the action selection problem. However,
		if the costs are non-convex in the capacity decisions then we may need to solve a mixed-integer nonlinear programming problem
		to update the first-stage decision in each iteration.
	\end{enumerate}
	
	Earlier in this paper, we compared the economic performance of NN-FVI
	to an inflexible two-stage MCIP model via out-of-sample tests. However,
	as indicated by \citet{zhao_decision_2018}, the complexity of the
	out-of-sample tests for ADP can be even greater than solving the original
	problem. One possible way to simplify the out-of-sample
	tests is to approximate the policy of the NN-FVI offline. For example,
	after solving the action selection problem for each state sample,
	we can do function fitting on these samples and then approximate the policy. When the policy is approximated offline, it does not add to
	the cost of the original NN-FVI algorithm. Then, we can use this approximate
	policy in the out-of-sample tests. In this case, the out-of-sample
	tests are more tractable and it is easier to implement the resulting
	policy in practice. Future work can consider how to approximate the
	policy to achieve high precision in out-of-sample tests.

	\section{Conclusion}
	
	In this paper, we solve a class of MDPs with a continuous state space and a large but finite and high-dimensional action space. An NN-FVI algorithm is proposed, where the value functions are approximated with a two-layer NN with ReLU activation functions. The consistency of NN-FVI is also formally proved. Since the action selection procedure of NN-FVI is time-consuming, we formulate it as a two-stage stochastic programming problem with a non-concave recourse function, and we design an MCD algorithm to solve it. We verify that the MCD algorithm converges to the global optimum in a finite number of iterations. We test our new algorithms on a capacity investment problem. Our numerical studies show that the proposed acceleration method significantly speeds up the action selection problem when compared with the brute-force method and the integer L-shape algorithm. 
    
    We solve discounted finite horizon MDPs in this paper, but our results on NN-FVI with MCD also extend to infinite horizon MDPs with minor modification. In addition, acceleration techniques that work for Benders decomposition may also be included to speed up our MCD algorithm. Some other possible future research directions include using NNs with different activation functions or multi-layer NNs along with new valid cuts. We may also extend our method to solve other problems, e.g. MCIP with non-convex capacity expansion costs.

	\section{Acknowledgment}
	{\small{}This research is supported by the National Natural Science Foundation of China (No. 72001141), and by the National Research Foundation, Prime Minister's Office, Singapore under its Campus for Research Excellence and Technological Enterprise (CREATE) programme.}
	
	%
	%
	%
		\appendix
		\label{sec:appendix}
		\section{Proofs of Main Mathematical Results}
		
		\subsubsection*{Proof of Proposition \ref{prop:Lip_Value_Function}}
		
		We first show that there exists $v_{\max}<\infty$ such that $|V_{t}\left(x\right)| \leq v_{\max}$ for all $x \in \mathbb{X}$ and $t=1,\ldots,T$. According to Assumption \ref{assu:MDP_finite_and_discounted}, there exists $r_{\max}<\infty$ such that $r_{\max}\ge\left|r_{t}\left(x,a\right)\right|$ for all $x\in\mathbb{X}$, $a\in\mathbb{A}$, and $t\in\mathbb{T}$. For $t=T$, we have $\left|V_{T}\left(x\right)\right|=\left|\max_{a\in\mathbb{A}}r_{T}\left(x,a\right)\right|\le r_{\max}$ for all $x\in \mathbb{X}$. For $t=1,\ldots,T-1$, we may take $v_{\max}\triangleq\max\left\lbrace\left(\sum_{\tau=1}^{T}\gamma^{\tau-1}\right)r_{\max},L_r\right\rbrace$ so that
		\[
		\left|V_{t}\left(x\right)\right|\le\max_{a\in\mathbb{A}}\left|r_{t}\left(x,a\right)+\gamma\int_{x^{\prime}\in\mathbb{X}}V_{t+1}\left(x^{\prime}\right)dP\left(x^{\prime}\left|x,a\right.\right)\right|\le v_{\max},\ \ \forall x\in\mathbb{X},\, t=1,\ldots,T-1.
		\]
		
		Now we check Lipschitz continuity of $V_{t}\left(\cdot\right)$. For $t=T$, we have 
		\[
		\left|V_{T}\left(x\right)-V_{T}\left(x^{\prime}\right)\right|=\left|\max_{a\in\mathbb{A}}r_{T}\left(x,a\right)-\max_{a^{\prime}\in\mathbb{A}}r_{T}\left(x^{\prime},a^{\prime}\right)\right|\le\max_{a\in\mathbb{A}}\left|r_{T}\left(x,a\right)-r_{T}\left(x^{\prime},a\right)\right|\le v_{\text{max}}\left\Vert x-x^{\prime}\right\Vert .
		\]
		In general, for $t<T$, we have{\small{}
			\begin{alignat*}{1}
			\left|V_{t}\left(x\right)-V_{t}\left(x^{\prime}\right)\right|\le & \left|\max_{a\in\mathbb{A}}\left(r_{t}\left(x,a\right)+\gamma\mathbb{E}V_{t+1}\left(y|x,a\right)\right)-\max_{a^{\prime}\in\mathbb{A}}\left(r_{t}\left(x^{\prime},a^{\prime}\right)+\gamma\mathbb{E}V_{t+1}\left(y^{\prime}|x^{\prime},a^{\prime}\right)\right)\right| \\
			\le & \max_{a\in\mathbb{A}}\left|\left(r_{t}\left(x,a\right)-r_{t}\left(x^{\prime},a\right)\right)\right|+\max_{a\in\mathbb{A}}\left|\gamma\mathbb{E}V_{t+1}\left(y\left|x,a\right.\right)-\gamma\mathbb{E}V_{t+1}\left(y^{\prime}\left|x^{\prime},a\right.\right)\right| \\
			\le & v_{\text{max}}\left\Vert x-x^{\prime}\right\Vert +\gamma\left(\max_{a\in\mathbb{A}}\int\left|V_{t+1}\left(y\right)\left(P\left(y\left|x,a\right.\right)-P\left(y\left|x^{\prime},a\right.\right)\right)\right|dy\right)\\
			\le & v_{\text{max}}\left\Vert x-x^{\prime}\right\Vert +\gamma v_{\text{max}}\max_{a\in\mathbb{A}}\int\left|P\left(y\left|x,a\right.\right)-P\left(y\left|x^{\prime},a\right.\right)\right|dy\\
			\le & \left(1+\gamma L_{P}\right)v_{\text{max}}\left\Vert x-x^{\prime}\right\Vert ,
			\end{alignat*}
		}using the fact that $V_{t}\left(\cdot\right)$ for $t\in\mathbb{T}$
		are bounded according to Assumption \ref{assu:Assumption2_Lip}. Thus,
		there exists $L_{v}=\left(1+\gamma L_{P}\right)v_{\text{max}}$
		such that $V_{t}\left(\cdot\right)$ are Lipschitz functions. \qed

		\subsubsection*{Proof of Lemma \ref{lem:Universal_Appro}}
		
		We only provide a sketch of proof here, and refer the interested readers
		to \citep{hornik_multilayer_1989} or \citep{sonoda_neural_2015}
		for the details. According to \citep[Corollary 2.2]{hornik_multilayer_1989},
		a two-layer NN is universal if the activation function
		of the network (denoted as $\Psi(x)$) is a squashing function. A
		squashing function satisfies the following properties according
		to \citep[Definition 2.3]{hornik_multilayer_1989}: $\Psi(x)$ is
		non-decreasing, $\lim_{x\rightarrow\infty}\Psi(x)=1$, and $\lim_{x\rightarrow-\infty}\Psi(x)=0$.
		
		ReLU does not satisfy the second property of a squashing function,
		but we can construct an equivalent network with a squashing activation
		function by taking a linear combination of two ReLUs. Suppose we have an
		ReLU network $\Gamma(x)$ with sufficient neurons. We may assume the
		number of neurons $J$ is even without loss of generality. For neuron
		$j\in\{1,\ldots,J/2\}$, we pick $j^{\prime}=J/2+j$ and specify
		its adjustable weights so that $u_{j^{\prime}}=u_{j}$, $u_{0j^{\prime}}=u_{0j}-1$,
		and $w_{j^{\prime}}=-w_{j}$. Then, we can construct an activation function by taking a linear combination of two neurons: 
		\[
		\hat{\Psi}_{j}(x)=\max\left\{ 0,u_{j}^{\top}x+u_{0j}\right\} -\max\left\{ 0,u_{j}^{\top}x+u_{0j}-1\right\} .
		\]
		Given $\hat{\Psi}_{j}(x)$, a new network $\hat{\Gamma}(x)$ with
		$J/2$ neurons is thus constructed. It follows that $\hat{\Gamma}$ is
		universal since $\hat{\Psi}_{j}(x)$ are squashing functions: $\hat{\Psi}(x)$
		is non-decreasing, $\lim_{x\rightarrow\infty}\hat{\Psi}(x)=1$, and
		$\lim_{x\rightarrow-\infty}\hat{\Psi}(x)=0$. Since the solutions
		of the weights of $\hat{\Gamma}$ are a subset of those of $\Gamma$,
		there always exists a network $\Gamma$ that is equivalent to $\hat{\Gamma}$.
		Thus, the two-layer network with ReLU is universal. \qed
		
		\subsubsection*{\label{sec:Proof-of-Consistency}Proof of Theorem \ref{prop:Consistency}}
		
		This result is derived from \citep[Corollary, 4]{munos_finite-time_2008}.
		Since we have shown that the two-layer NN is a universal
		approximator in Lemma \ref{lem:Universal_Appro}, the only thing we
		need to prove now is that the MDP presented in this paper satisfies
		the assumptions of MDP regularity and uniformly stochastic transitions
		from \citep[Corollary, 4]{munos_finite-time_2008}. Let $P_{R}\left(\cdot|x,a\right)$
		be the probability distribution of reward $r_{t}$ given state $x\in\mathbb{X}$
		and action $a\in\mathbb{A}$ in time $t$. The assumptions stated in \citep{munos_finite-time_2008}
		are presented below.
		\begin{assumption}
		    \label{assu:mdp_regularity}
			{[}MDP regularity{]} The MDP satisfies the following conditions: $\mathbb{X}$
			is a bounded, closed subset of some Euclidean space, $\mathbb{A}$ is
			finite and the discount factor $\gamma$ satisfies $0<\gamma<1$.
			The reward kernel is such that the immediate reward function is a
			bounded measurable function with bound $r_{\max}$. Further, the support
			of $P_{R}\left(\cdot|x,a\right)$ is included in $\left[-r_{\max},r_{\max}\right]$
			independently of $\left(x,a\right)\in\mathbb{X}\times\mathbb{A}$.
		\end{assumption}
		
		\begin{assumption}
		    \label{assu:uniform_transitions}
			{[}Uniformly stochastic transitions{]} For all $x\in\mathbb{X}$ and
			$a\in\mathbb{A}$, $P\left(\cdot|x,a\right)$ is absolutely
			continuous w.r.t. $\mu$. Also, the Radon-Nikodym derivative of $P$ w.r.t.
			$\mu$ is bounded uniformly with bound $C_{\mu}:$
			\[
			C_{\mu}\triangleq\sup_{x\in\mathbb{X},a\in\mathbb{A}}\left\Vert \frac{dP\left(\cdot|x,a\right)}{d\mu}\right\Vert _{\infty}<+\infty.
			\]
		\end{assumption}
		
		We first show that the MDP presented in Section \ref{sec:Preliminaries}
		satisfies Assumption \ref{assu:mdp_regularity} According to Assumption \ref{assu:MDP_finite_and_discounted},
		there exists $r_{\max}$ such that $r_{t}\in\left[-r_{\max},r_{\max}\right]$.
		As $r_{t}\left(x,a\right)$ is deterministic and bounded, the support of $P_{R}\left(\cdot|x,a\right)$ is included in the bounded
		set $\left[-r_{\max},r_{\max}\right]$ which is independent of $\left(x,a\right)\in\mathbb{X}\times\mathbb{A}$.
		
		Now we show that Assumption \ref{assu:uniform_transitions} is satisfied. First, according to
		\citep{munos_finite-time_2008}, Assumption \ref{assu:uniform_transitions} is equivalent to assuming
		that the transition kernel admits a uniformly bounded density when
		$\mu$ is the Lebesgue measure over $\mathbb{X}$, which is essentially
		Assumption \ref{assu:MDP_lip_transition}. Therefore, all of the required conditions
		are satisfied. \qed
		
		\subsubsection*{Proof of Proposition \ref{prop:valid_inquality_neg}}
		
		For all $j\in\mathcal{J}_{t+1}^{-}$, we
		have $w_{j\left(t+1\right)}\le0$. For all $s\in\mathcal{S}_{2}$,
		denote 
		\[\eta_{j}^{s}=w_{j\left(t+1\right)}\max\left\{ u_{j\left(t+1\right)}^{\top}\left(A\left(x_{t},\xi_{t}^{s}\right)+B\left(x_{t},\xi_{t}^{s}\right)a_{t}\right)+u_{0j\left(t+1\right)},0\right\}. \]
		Then, we have 
		\begin{alignat*}{1}
		\eta_{j}^{s} & =\min\left\{ w_{j\left(t+1\right)}\left(u_{j\left(t+1\right)}^{\top}\left(A\left(x_{t},\xi_{t}^{s}\right)+B\left(x_{t},\xi_{t}^{s}\right)a_{t}\right)+u_{0j\left(t+1\right)}\right),0\right\} \\
		& \le \min\left\lbrace\left(\phi_{j}^{ms}\right)^{\top}a_{t}+\phi_{0j}^{ms}, 0\right\rbrace\\
		& \le\left(\phi_{j}^{ms}\right)^{\top}a_{t}+\phi_{0j}^{ms},\ \ \forall a_{t}\in\mathbb{A}.
		\end{alignat*}
		The first equality holds because $w_{j\left(t+1\right)}\le0$. The second line follows from the definitions
		of $\phi_{nj}^{ms}$ and $\phi_{0j}^{ms}$ given a specific $a_{t}^{m}\in\mathbb{A}$: if $u_{j\left(t+1\right)}^{\top}f\left(x_t,a_t^m,\xi_t^s\right)+u_{0j\left(t+1\right)}>0$,
		the inequality holds trivially; if $u_{j\left(t+1\right)}^{\top}f\left(x_t,a_t^m,\xi_t^s\right)+u_{0j\left(t+1\right)}\le0$,
		then $\phi_{nj}^{ms}$ and $\phi_{0j}^{ms}$ are all zero and so the	inequality holds because the first line is always non-positive. Summing
		up over all neurons $j\in\mathcal{J}_{t+1}^{-}$ and then taking expectations
		on both sides of the inequality, we have 
		\begin{equation*}
		\eta^{-} \le
		\frac{1}{S_2}\sum_{s\in\mathcal{S}_{2}}\sum_{j\in\mathcal{J}_{t+1}^{-}}\eta_{j}^{s}
		\le
		\frac{1}{S_2}\sum_{s\in\mathcal{S}_{2}}\sum_{j\in\mathcal{J}_{t+1}^{-}}\left(\left(\phi_{j}^{ms}\right)^{\top}a_{t}+\phi_{0j}^{ms}\right),\ \ \forall a_{t}\in\mathbb{A}.
		\end{equation*}
		Thus, according to the definition of the hypograph, we have
		\[
		\text{hypo}_{t}^{-}\subset\left\{ \left(a_{t},\eta\right)\in\mathbb{R}^{N_2+1}\left|
		\eta^{-}\le\frac{1}{S_2}\sum_{s\in\mathcal{S}_{2}}\sum_{j\in\mathcal{J}_{t+1}^{-}}\left(\left(\phi_{j}^{ms}\right)^{\top}a_{t}+\phi_{0j}^{ms}\right)\right.\right\} ,
		\]
		and so Eq. \eqref{eq:validcut_neg} is a valid inequality. \qed
		
		\subsubsection*{Proof of Proposition \ref{prop:valid_inquality_pos}}
		To prove that Eq. \eqref{eq:validcut_pos} is a valid inequality for $\text{hypo}_{t}^{+}$,
		we need to prove 
		\begin{alignat*}{1}
		\sum_{n=1}^{N_{2}}\theta_{nj}^{s}a_{nt}+\theta_{0j}^{s} & \ge w_{j\left(t+1\right)}\max\left\{ u_{j\left(t+1\right)}^{\top}f\left(x_{t},a_{t},\xi_{t}^{s}\right)+u_{0j\left(t+1\right)},0\right\} ,\\
		& =w_{j\left(t+1\right)}\max\left\{ \sum_{n=1}^{N_{2}}\Gamma_{njs}^{\left(1\right)}a_{nt}+\Gamma_{js}^{\left(2\right)},0\right\} ,\ \ \forall a_{t}\in\mathbb{A}.
		\end{alignat*}
		
		(i) If $\sum_{n\in\mathbf{N}_{js}^{-}}\Gamma_{njs}^{\left(1\right)}\bar{a}_{n}+\Gamma_{js}^{\left(2\right)}>0$,
		then we have 
		\begin{equation*}
		\sum_{n=1}^{N_{2}}\theta_{nj}^{s}a_{nt}+\theta_{0j}^{s} =w_{j(t+1)} \left(\sum_{n=1}^{N_{2}}\Gamma_{njs}^{\left(1\right)}a_{nt}+\Gamma_{js}^{\left(2\right)}\right)
		= w_{j\left(t+1\right)}\max\left\{ \sum_{n=1}^{N_{2}}\Gamma_{njs}^{\left(1\right)}a_{nt}+\Gamma_{js}^{\left(2\right)},0\right\}.
		\end{equation*}
		The first equality is derived from the definitions of $\theta_{nj}^{s}$
		and $\theta_{0j}^{s}$. The second equality holds because $\sum_{n=1}^{N_{2}}\Gamma_{njs}^{\left(1\right)}a_{nt}+\Gamma_{js}^{\left(2\right)} \ge \sum_{n\in\mathbf{N}_{js}^{-}}\Gamma_{njs}^{\left(1\right)}\bar{a}_{n}+\Gamma_{js}^{\left(2\right)} > 0$ for all $a_t\in\mathbb{A}$. 
		
		(ii) If $\sum_{n\in\mathbf{N}_{js}^{+}}\Gamma_{njs}^{\left(1\right)}\bar{a}_{n}+\Gamma_{js}^{\left(2\right)}<0$, the result holds trivially	because $\sum_{n\in\mathbf{N}_{js}^{+}}\Gamma_{njs}^{\left(1\right)}\bar{a}_{n}+\Gamma_{js}^{\left(2\right)}$
		is the maximum value of $\sum_{n=1}^{N_{2}}\Gamma_{njs}^{\left(1\right)}a_{nt}+\Gamma_{js}^{\left(2\right)}$,
		such that $\max\left\{ \sum_{n=1}^{N_{2}}\Gamma_{njs}^{\left(1\right)}a_{nt}+\Gamma_{js}^{\left(2\right)},0\right\} =0$.
		
		(iii) We prove that the inequality is valid when $\sum_{n\in\mathbf{N}_{js}^{-}}\Gamma_{njs}^{\left(1\right)}\bar{a}_{n}+\Gamma_{js}^{\left(2\right)}\le0$
		and $\sum_{n\in\mathbf{N}_{js}^{+}}\Gamma_{njs}^{\left(1\right)}\bar{a}_{n}+\Gamma_{js}^{\left(2\right)}\ge0$.
		First, since $0 \le a_{nt}\le\bar{a}_{n}$ for all $n=1,\ldots,N_{2}$, we
		have 
		\begin{equation*}
		\sum_{n=1}^{N_{2}}\Gamma_{njs}^{\left(1\right)}a_{nt}-\sum_{n\in\mathbf{N}_{js}^{+}}\Gamma_{njs}^{\left(1\right)}\bar{a}_{n}=\sum_{n\in\mathbf{N}_{js}^{+}}\Gamma_{njs}^{\left(1\right)}\left(a_{nt}-\bar{a}_{n}\right)+\sum_{n\in\mathbf{N}_{js}^{-}}\Gamma_{njs}^{\left(1\right)}a_{nt}\le0,
		\end{equation*}
		according to the definitions of $\mathbf{N}_{js}^{+}$ and $\mathbf{N}_{js}^{-}$.
		Then, we have 
		\begin{flalign*}
		\vartheta_{j}^{s} & =\frac{\sum_{n\in\mathbf{N}_{js}^{+}}\Gamma_{njs}^{\left(1\right)}\bar{a}_{n}+\Gamma_{js}^{\left(2\right)}}{\sum_{n=1}^{N_{2}}\left|\Gamma_{njs}^{\left(1\right)}\right|\bar{a}_{n}},\nonumber \\
		& =\frac{\sum_{n\in\mathbf{N}_{js}^{+}}\Gamma_{njs}^{\left(1\right)}\bar{a}_{n}-\sum_{n\in\mathbf{N}_{js}^{-}}\Gamma_{njs}^{\left(1\right)}\bar{a}_{n}+\sum_{n\in\mathbf{N}_{js}^{-}}\Gamma_{njs}^{\left(1\right)}\bar{a}_{n}+\Gamma_{js}^{\left(2\right)}}{\sum_{n\in\mathbf{N}_{js}^{+}}\Gamma_{njs}^{\left(1\right)}\bar{a}_{n}-\sum_{n\in\mathbf{N}_{js}^{-}}\Gamma_{njs}^{\left(1\right)}\bar{a}_{n}}\nonumber \\
		& =1+\frac{\sum_{n\in\mathbf{N}_{js}^{-}}\Gamma_{njs}^{\left(1\right)}\bar{a}_{n}+\Gamma_{js}^{\left(2\right)}}{\sum_{n\in\mathbf{N}_{js}^{+}}\Gamma_{njs}^{\left(1\right)}\bar{a}_{n}-\sum_{n\in\mathbf{N}_{js}^{-}}\Gamma_{njs}^{\left(1\right)}\bar{a}_{n}}\nonumber \\
		& \le1,
		\end{flalign*}
		where the inequality holds because $\sum_{n\in\mathbf{N}_{js}^{-}}\Gamma_{njs}^{\left(1\right)}\bar{a}_{n}+\Gamma_{js}^{\left(2\right)}\le0$. Note that $w_{j\left(t+1\right)} >0$ for all $j \in \mathcal{J}_{t+1}^{+}$ by definition. Therefore, for any $a_{t}\in\mathbb{A}$ we have	
		\begin{align*}
		& \left(\sum_{n=1}^{N_{2}}\theta_{nj}^{s}a_{nt}+\theta_{0j}^{s}\right)\big/w_{j\left(t+1\right)}\\
		= & \vartheta_{j}^{s}\left(\sum_{n=1}^{N_{2}}\Gamma_{njs}^{\left(1\right)}a_{nt}-\sum_{n\in\mathbf{N}_{js}^{-}}\Gamma_{njs}^{\left(1\right)}\bar{a}_{n}\right)\\
		= & \vartheta_{j}^{s}\left(\sum_{n=1}^{N_{2}}\Gamma_{njs}^{\left(1\right)}a_{nt}-\sum_{n\in\mathbf{N}_{js}^{+}}\Gamma_{njs}^{\left(1\right)}\bar{a}_{n}+\sum_{n=1}^{N_{2}}\left|\Gamma_{njs}^{\left(1\right)}\right|\bar{a}_{n}\right)\\
		= & \vartheta_{j}^{s}\left(\sum_{n=1}^{N_{2}}\Gamma_{njs}^{\left(1\right)}a_{nt}-\sum_{n\in\mathbf{N}_{js}^{+}}\Gamma_{njs}^{\left(1\right)}\bar{a}_{n}\right)+\sum_{n\in\mathbf{N}_{js}^{+}}\Gamma_{njs}^{\left(1\right)}\bar{a}_{n}+\Gamma_{js}^{\left(2\right)}\\
		= & \vartheta_{j}^{s}\left(\sum_{n=1}^{N_{2}}\Gamma_{njs}^{\left(1\right)}a_{nt}-\sum_{n\in\mathbf{N}_{js}^{+}}\Gamma_{njs}^{\left(1\right)}\bar{a}_{n}\right)-\left(\sum_{n=1}^{N_{2}}\Gamma_{njs}^{\left(1\right)}a_{nt}-\sum_{n\in\mathbf{N}_{js}^{+}}\Gamma_{njs}^{\left(1\right)}\bar{a}_{n}\right)+\sum_{n=1}^{N_{2}}\Gamma_{njs}^{\left(1\right)}a_{nt}+\Gamma_{js}^{\left(2\right)}\\
		= & \left(\vartheta_{j}^{s}-1\right)\left(\sum_{n=1}^{N_{2}}\Gamma_{njs}^{\left(1\right)}a_{nt}-\sum_{n\in\mathbf{N}_{js}^{+}}\Gamma_{njs}^{\left(1\right)}\bar{a}_{n}\right)+\sum_{n=1}^{N_{2}}\Gamma_{njs}^{\left(1\right)}a_{nt}+\Gamma_{js}^{\left(2\right)}\\
		\ge & \sum_{n=1}^{N_{2}}\Gamma_{njs}^{\left(1\right)}a_{nt}+\Gamma_{js}^{\left(2\right)}.
		\end{align*}
		The inequality holds because we have $\vartheta_{j}^{s}\le1$ and $\sum_{n=1}^{N_{2}}\Gamma_{njs}^{\left(1\right)}a_{nt}-\sum_{n\in\mathbf{N}_{js}^{+}}\Gamma_{njs}^{\left(1\right)}\bar{a}_{n}\le0$ when $\sum_{n\in\mathbf{N}_{js}^{-}}\Gamma_{njs}^{\left(1\right)}\bar{a}_{n}+\Gamma_{js}^{\left(2\right)} \le 0$ for all $j\in\mathcal{J^{+}}$. 
		
		The last step is to show that $\sum_{n=1}^{N_{2}}\theta_{nj}^{s}a_{nt}+\theta_{0j}^{s}\ge0$
		for all $a_{t}\in\mathbb{A}$:
		\begin{alignat*}{1}
		\sum_{n=1}^{N_{2}}\theta_{nj}^{s}a_{nt}+\theta_{0j}^{s}= & w_{j(t+1)} \vartheta_j^s \left(\sum_{n=1}^{N_{2}}\Gamma_{njs}^{\left(1\right)}a_{nt}-\sum_{n\in\mathbf{N}_{js}^{-}}\Gamma_{njs}^{\left(1\right)}\bar{a}_{n}\right)\\
		= & w_{j(t+1)} \vartheta_j^s \left(\sum_{n\in\mathbf{N}_{js}^{+}}\Gamma_{njs}^{\left(1\right)}a_{nt}+\sum_{n\in\mathbf{N}_{js}^{-}}\Gamma_{njs}^{\left(1\right)}\left(a_{nt}-\bar{a}_{n}\right)\right) \ge 0.
		\end{alignat*}
		The inequality holds because $\vartheta_j^s \ge 0$ when $\sum_{n\in\mathbf{N}_{js}^{+}}\Gamma_{njs}^{\left(1\right)}\bar{a}_{n}+\Gamma_{js}^{\left(2\right)}\ge0$, $\Gamma_{njs}^{\left(1\right)}a_{nt}\ge0$ for all $n\in\mathbf{N}_{js}^{+}$, and $\Gamma_{njs}^{\left(1\right)}\left(a_{nt}-\bar{a}_{n}\right)\ge0$ for all $n\in\mathbf{N}_{js}^{-}$.	\qed
		
		\subsubsection*{Proof of Theorem \ref{prop:finite-conveg}}
		
		According to \citep[Proposition 2]{laporte_integer_1993} and Corollary \ref{coro:valid_cut_all}, Cuts \eqref{eq:LLCuts_MP} and \eqref{eq:validcut_all2}
		are all valid. Since the action space $\mathbb{A}$ is finite, Problem
		\eqref{eq:FP} has finitely many feasible solutions. As a result,
		there are only finitely many cuts that can be added to Problem \eqref{eq:FP}. Thus, the algorithm will	converge to the global optimal solution in a finite number of iterations.
		\qed 
		
		\subsubsection*{Proof of Proposition \ref{prop:check_assumptions}}
		
		We need to check Assumptions \ref{assu:MDP_finite_and_discounted} and
		\ref{assu:MDP_linear_transition}. Assumption \ref{assu:MDP_lip_reward} holds trivially as the unit revenue and unit penalty cost in Eq. \eqref{eq:Lp_obj} are finite. Assumption \ref{assu:exogenous} holds according to \citep[Proposition 8.6]{kallenberg_foundations_2002}. Assumption \ref{assu:MDP_lip_transition}
		holds trivially if Assumption \ref{assu:MDP_linear_transition} is
		satisfied. 
		
		To verify Assumption \ref{assu:MDP_finite_and_discounted}, we need to confirm that
		the reward function given by Eq. \eqref{eq:cost_function} is bounded.
		As $K_{t}\le K^{\max}$ and $d_{t}\le D^{\max}$ for all $t\in\mathbb{T}$,
		an upper and lower bound on $Q_{t}\left(K_{t-1},d_{t}\right)$ follows
		if we have unlimited or no capacity: 
		\[
		-\sum_{i\in\mathcal{I}}\max_{t\in\mathbb{T}}b_{it}D_{i}^{\max}\le Q_{t}\left(K_{t-1},d_{t}\right)\le\sum_{i\in\mathcal{I}}\max_{t\in\mathbb{T},n\in\mathcal{N}}\hat{r}_{int}D_{i}^{\max}.
		\]
		Therefore, 
		\[
		\left|Q_{t}\left(K_{t-1},d_{t}\right)\right|\le\sum_{i\in\mathcal{I}}\max_{t\in\mathbb{T},n\in\mathcal{N}}\left(\hat{r}_{int}+b_{it}\right)D_{i}^{\max}.
		\]
		Similarly, upper bounds on the capacity adjustment cost
		follow by assuming that the capacity is changed from zero to $K^{\max}$
		or the reverse: 
		\[
		-\sum_{n\in\mathcal{N}}\max_{t\in\mathbb{T}}q_{nt}^{-}K_{n}^{\max}\le\sum_{n\in\mathcal{N}}\max\left\{ -q_{nt}^{-}\left(K_{n\left(t-1\right)}-K_{nt}\right),q_{nt}^{+}\left(K_{nt}-K_{n\left(t-1\right)}\right)\right\} \le\sum_{n\in\mathcal{N}}\max_{t\in\mathbb{T}}q_{nt}^{+}K_{n}^{\max}.
		\]
		Since we have assumed $q_{nt}^{-}\le q_{nt}^{+}$, it follows that
		\[
		\left|\sum_{n\in\mathcal{N}}\max\left\{ -q_{nt}^{-}\left(K_{n\left(t-1\right)}-K_{nt}\right),q_{nt}^{+}\left(K_{nt}-K_{n\left(t-1\right)}\right)\right\} \right|\le\sum_{n\in\mathcal{N}}\max_{t\in\mathbb{T}}q_{nt}^{+}K_{n}^{\max}.
		\]
		Then, for all $K_{t}\in\mathbb{K}$ and $t\in\mathbb{T}$, we have 
		\begin{align*}
		\left|r_{t}\left(\left(K_{t-1},d_{t}\right),K_{t}\right)\right| & \le\left|Q_{t}\left(K_{t-1},d_{t}\right)\right|+\left|\sum_{n\in\mathcal{N}}\max\left\{ -q_{nt}^{-}\left(K_{n\left(t-1\right)}-K_{nt}\right),q_{nt}^{+}\left(K_{nt}-K_{n\left(t-1\right)}\right)\right\} \right|\\
		& \le\sum_{n\in\mathcal{N}}\max_{t\in\mathbb{T}}q_{nt}^{+}K_{n}^{\max}+\sum_{i\in\mathcal{I}}\max_{t\in\mathbb{T},n\in\mathcal{N}}\left(\hat{r}_{int}+b_{it}\right)D_{i}^{\max}\\
		& <\infty.
		\end{align*}
		
		For Assumption \ref{assu:MDP_linear_transition}, we let $0_{I\times N}$ be the $\left(I\times N\right)$-dimensional matrix of zeros and $\text{diag}\left(1_N\right)$ be the $(N\times N)$-dimensional identity matrix. 		According to \citep[Proposition 8.6]{kallenberg_foundations_2002}, since $\{D_t\}_{t\in\mathbb{T}}$ is Markov there exist measurable functions $f_t^D$ and i.i.d. $Uniform(0,1)$ random variable $\xi_t$ such that $d_{t} = f_{t}^D \left(d_{t-1}, \xi_{t}\right)$. If we let $x_t \triangleq \left(K_{t-1}, d_t\right)$ be the state variable, then the transition function follows by
		\[
		\left[\begin{array}{c}
		K_{t}\\
		d_{t+1}
		\end{array}\right]=\left[\begin{array}{c}
		K_{t-1}\\
		f_{t+1}^D\left(d_{t}, \xi_t\right)
		\end{array}\right]+\left[\begin{array}{c}
		\text{diag}\left(1_{N}\right)\\
		0_{I\times N}
		\end{array}\right]\Delta K_{t}, \ \ \forall t \in \mathbb{T}\backslash\{T\},
		\]
		which is a linear function with respect to $\Delta K_t$. \qed
		
		\subsubsection*{Proof of Proposition \ref{prop:equivalence_value_functions}}
		
		According to the definition of $\bar{V}_{T}\left(\cdot\right)$, we
		have $\bar{V}_{T}\left(\cdot\right)=V_{T}\left(\cdot\right)$ at all arguments in $\left(K_{T-1},d_{T}\right)\in\mathbb{X}$. Suppose we have
		$\bar{V}_{t+1}\left(\cdot\right)=V_{t+1}\left(\cdot\right)$ for all
		$\left(K_{t},d_{t+1}\right)\in\mathbb{X}$. Since the demand transitions
		are independent of $K_{t-1}$ and $K_{t}$, we have $\mathbb{E}\left[\left.\bar{V}_{t+1}\left(K_{t},D_{t+1}\right)\right|d_{t}\right]=\mathbb{E}\left[\left.V_{t+1}\left(K_{t},D_{t+1}\right)\right|d_{t}\right]$
		for all $K_{t}\in\mathbb{K}$. Therefore, $\bar{V}_{t}\left(K_{t-1},d_{t}\right)=V_{t}\left(K_{t-1},d_{t}\right)$
		for all $\left(K_{t-1},d_{t}\right)\in\mathbb{X}$ according to Eq.
		\eqref{eq:ValueFunction} and Eq. \eqref{eq:ValueFunction-Appro}.
		The desired result follows by backward induction. \qed
		
		\subsubsection*{Proof of Proposition \ref{prop:Piecewise Value Function}}
		
		First, we know that the property of piecewise linearity is preserved under finite
		summation and the max/min of a finite collection. Now, observe that
		$Q_{t}\left(K_{t-1},d_{t}\right)$ is a piecewise linear function
		when defined over $\left(K_{t-1},d_{t}\right)\in\bar{\mathbb{X}}$.
		This observation follows by transforming $Q_{t}\left(K_{t-1},d_{t}\right)$
		into its dual problem. Since $Q_{t}\left(K_{t-1},d_{t}\right)$ is a linear
		programming problem and the optimal value is finite and attained for
		all $\left(K_{t-1},d_{t}\right)\in\bar{\mathbb{X}}$, strong duality
		holds (see e.g. \citep[Chapter 4]{bertsimas1997introduction}) and
		we have:
		\begin{alignat*}{1}
		Q_{t}\left(K_{t-1},d_{t}\right)=\min_{\psi_{n},\lambda_{i}\ge0} & \sum\limits _{i\in\mathcal{I}}\left(b_{it}-\lambda_{i}\right)d_{it}-\sum_{n\in\mathcal{N}}\psi_{n}K_{n\left(t-1\right)}\\
		\text{s. t. } & \left(\psi_{n}+\lambda_{i}-\hat{r}_{int}-b_{it}\right)\ge0,\ \ \forall i\in\mathcal{I},n\in\mathcal{N}.
		\end{alignat*}
		Let $\Lambda$ denote the feasible set of $\left(\psi,\lambda\right)$
		for the above problem, where $\lambda=\left(\lambda_{i}\right)_{i\in\mathcal{I}}$
		and $\psi=\left(\psi_{n}\right)_{n\in\mathcal{N}}$. Also let $\text{ext}\left(\Lambda\right)$ be the set of extreme points of $\Lambda$. Since the dual problem is an LP, there is an optimal solution $\left(\psi^{*},\lambda^{*}\right)$
		among $\text{ext}\left(\Lambda\right)$ \citep[Chapter 3]{bertsimas1997introduction},
		and so the above problem is equivalent to 
		\[
		\min_{\left(\psi,\lambda\right)\in\text{ext}\left(\Lambda\right)}\sum\limits _{i\in\mathcal{I}}\left(b_{it}-\lambda_{i}\right)d_{it}-\sum_{n\in\mathcal{N}}\mu_{n}K_{n\left(t-1\right)}.
		\]
		It follows that $Q_{t}\left(K_{t-1},d_{t}\right)$ is piecewise linear
		and concave in $K_{t-1}\in\bar{\mathbb{K}}$ since it is the min of
		a finite collection of linear functions. Since $\bar{V}_{T}\left(K_{T-1},d_{T}\right)=r_{T}\left(\left(K_{T-1},d_{T}\right),\boldsymbol{0}_{N}\right)$
		and $r_{T}\left(\cdot\right)$ is concave in $K_{T-1}$, $\bar{V}_{T}\left(K_{T-1},d_{T}\right)$
		is piecewise linear for all $\left(K_{T-1},d_{T}\right)\in\bar{\mathbb{X}}$
		and concave in $K_{T-1}$.
		
		To prove (ii), for the induction step suppose that $\bar{V}_{t+1}\left(K_{t},D_{t+1}\right)$ is piecewise linear. We have
		\[
		\bar{V}_{t}\left(K_{t-1},d_{t}\right)=\max\limits _{K_{t}\in\mathbb{K}}\left\{ r_t\left(\left(K_{t-1},d_{t}\right),K_{t}\right)+\gamma\mathbb{E}\left[\bar{V}_{t+1}\left(K_{t},D_{t+1}\right)\left|d_{t}\right.\right]\right\} ,\ \ \forall t\in\mathbb{T}\backslash\left\{ T\right\} .
		\]
		Since $D_{t+1}\in\mathbb{D}$ is finite, $\mathbb{E}\left[\bar{V}_{t+1}\left(K_{t},D_{t+1}\right)\left|d_{t}\right.\right]$
		is piecewise linear as it is a finite sum of piecewise linear functions.
		Therefore, $\bar{V}_{t}\left(K_{t-1},d_{t}\right)$ is piecewise linear
		as the max of a finite set of piecewise linear functions, and the desired result follows. \qed
		
		\section{Formulation of Integer Optimality Cuts}
		
		To formulate the integer optimality cuts, we first define a set of indices $\mathcal{L}_n \triangleq \{0,1,\ldots,L_n\}$ for all $n=1,\ldots,N_2$, such that $2^{L_n-1} \le \bar{a}_n \le 2^{L_n}$ for all $n=1,\ldots,N_2$. Then, we define binary variables $\alpha_{nlt}\in\{0,1\}$ for all $l\in\mathcal{L}_n$ such that $a_{nt} = \sum_{l\in\mathcal{L}_n} 2^l \alpha_{nlt}$. We let $\alpha_{nlt}^m$ denote the binary variables corresponding to the optimal action $a_{nt}^m$ in the $m$th iteration of MCD. We also let
		$$
		\mathcal{Z}_{1}^{m} \triangleq \left\lbrace (n, l) \left| \alpha_{nlt}^m=1, n\in\{1,\ldots,N_2\},\, l\in\mathcal{L}_n \right.\right\rbrace
		$$
		and
		$$\mathcal{Z}_{0}^{m} \triangleq \left\lbrace (n, l) \left| \alpha_{nlt}^m=0, n\in\{1,\ldots,N_2\},\, l\in\mathcal{L}_n \right.\right\rbrace
		$$
		be the index sets for the binary variables that are equal to one and zero, respectively. Then, $\zeta^{m}\left(a_{t}\right)$ can be expressed as:
		\[
		\zeta^m\left(a_t\right) = \left| \mathcal{Z}^m_1 \right| - \left[\sum_{(n,l)\in\mathcal{Z}^m_1}\alpha_{nlt} - \sum_{(n,l)\in\mathcal{Z}^m_0}\alpha_{nlt} \right],
		\]
		where $\left| \mathcal{Z}^m_1 \right|$ is the cardinality of the set $\mathcal{Z}^m_1$. If $a_t=a_t^m$, then $\zeta^m\left(a_t\right)=0$; otherwise $\zeta^m\left(a_t\right)\ge1$. Other possible formulations for $\zeta^m\left(a_t\right)$ can also be found in \citep{laporte_integer_1993}.
	
	
	\bibliographystyle{informs2014} 
	\bibliography{MyLibrary} 
	
	
\end{document}